\documentclass{article}

\pdfoutput=1
\usepackage{arxiv}

\usepackage{algorithmicx}

\usepackage{algorithm}
\usepackage{algpseudocode}
\usepackage[caption=false]{subfig}
\usepackage{mathtools}
\usepackage{array}    
\usepackage{rotating}

\usepackage[utf8]{inputenc} 
\usepackage[T1]{fontenc}    
\usepackage{hyperref}       
\usepackage{url}            
\usepackage{booktabs}       
\usepackage{amsfonts}       
\usepackage{nicefrac}       
\usepackage{microtype}      
\usepackage{lipsum}		

\usepackage{amsopn}

\ifpdf
  \DeclareGraphicsExtensions{.eps,.pdf,.png,.jpg}
\else
  \DeclareGraphicsExtensions{.eps}
\fi

\newcommand{\figref}[1]{Figure~#1}

\newcommand{\Algref}[1]{Algorithm~#1}
\newcommand{\sectref}[1]{Section~#1}
\newcommand{\subsectref}[1]{Subsection~#1}
\newcommand{\equaref}[1]{#1}
\newcommand{\Equaref}[1]{Equation~#1}

\title{LOCAL TIME STEPPING METHODS AND DISCONTINUOUS GALERKIN METHODS APPLIED TO DIFFUSION ADVECTION REACTION EQUATIONS}

\author{
  Assionvi H. Kouevi\thanks{Department of Applied Mathematics and computer sciences, Heriot-Watt University
    .} \\
  Dept.  of Applied Maths and computer sciences\\
  Heriot-Watt University\\
  Edinburgh, EH14 4AS \\
  \url{https://www.hw.ac.uk/} \\
   \And
 Gabriel J. Lord\footnotemark[1] \\
  Dept.  of Applied Maths and computer sciences\\
  Heriot-Watt University\\
  Edinburgh, EH14 4AS \\
  \url{www.macs.hw.ac.uk/\~gabriel/} \\
}

\begin{document}
\maketitle

\begin{abstract}
This paper is focussed on the numerical resolution of diffusion
advection and reaction equations (DAREs) with special features (such
as fractures, walls, corners, obstacles or point loads) which
globally, as well as locally, have important effects on the
solution. We introduce a multilevel and local time solver of DAREs
based on the discontinuous Galerkin (DG) method for the spatial
discreization and time stepping methods
such as exponential time differencing (ETD), exponential Rosenbrock
(EXPR) and implicit Euler (Impl) methods. The efficiency of our solvers
is shown with several experiments on cyclic voltammetry models and
fluid flows through domains with fractures.
\end{abstract}

\keywords{Local time stepping methods \and discontinuous Galerkin method \and implicit
Euler method \and exponential time differencing method \and exponential
Rosenbrock method \and diffusion advection and reaction \and cyclic voltammetry \and domain with
fracture \and unstructured mesh}

\section{Introduction}
We examine numerical methods for diffusion advection and reaction
equations of the form  
\begin{equation}\label{DARE}
  \dfrac{\partial C}{\partial t}  - \nabla \cdot (D\nabla C)+
  \nabla\cdot\textbf{v}C = R(C),  \text{ given } C(x,0)=C_0(x),
  \quad x \in \Omega,
\end{equation}
where, for example, $C$ is a concentration of a solute, where $R(C)$ is a reaction term, $D$ is a diffusivity and $\textbf{v}$ is a given velocity. Often, there are special features (e.g. fractures, walls, corners, obstacles, electrodes, point loads or irregular material interfaces), which affect locally the flow and transport of the solute. To accurately capture such local behaviour numerically, spatial local refinement is necessary. However, this requires a  reduction of the time step, $\Delta t$, for stability (while using the explicit time integrator) and for accuracy (while using both explicit and implicit time integrators). Unfortunately, when applied uniformly on all the simulation domain, $\Omega$, the reduced time step leads to an unacceptable large CPU time, making the use of local time stepping (LTS) methods highly desirable. The key feature of LTS methods is to split the solution domain $\Omega$ into several sub-domains $\Omega_i$ each with a time step $\Delta t_i$ as large as possible for efficiency. According to \cite{LTS1}, a LTS method is efficient if it ensures accuracy of the solution, i.e. the solution has to be more accurate than the one obtained with a global coarse mesh, and in addition leads to reduced CPU time compared to the one obtained when using a small time step on the whole domain. We give a short review of LTS methods in the context of the DG applximation.

The LTS methods have their roots in the work of Rice \cite{LTS2}, who in 1960 developped the so-called multirate Runge-Kutta methods for a two scale system of ordinary differential equations (ODEs). The multirate approach, for ODEs, was then combined with linear multistep integrators in 1984 by Gear et al. \cite{gear1984multirate} to improve the accuracy; and their stability properties were analyzed in 1989 by Skelboe et al \cite{skelboe1989stability}. These multirate
methods were based on the static partitioning of the domain $\Omega$ generated from the  priori knowledge of the physics of the problem. This limitation was overcomed in 1997 by Engstler et al. \cite{engstler1997m} when they introduced the multirate extrapolation methods for ODEs, based on Richardson extrapolation. Due to their unconditional stability, the (DG) method was also used to handle the time refinement of
ODEs.  The DG method, denoted DG($q$) while using the polynomials of degree $q\in \mathbb{N}$,  applied to ODEs was first studied in 1974 by LeSaint et al. \cite{Lesaint19}.  It was proven to be strongly A-stable of order $2q+1$ by the authors. Note that the case $q=0$ is equivalent to implicit Euler scheme. Adaptive error control was introduced in \cite{johnson1988error} and more recently considered in \cite{estep1995posteriori, bottcher1996adaptive}. For more insight on
the DG method applied to ODEs, see for example
\cite{cockburn2000DGreview}. 

In the case of partial differential  equations (PDEs), several schemes using the DG method have been developed to handle space and time refinement problems. The finite element (FE) method in space followed by the DG method in time was used in \cite{akrivis2004,chrysafinos2006, chrysafinos2006er, eriksson1991, John1, estep1993disc, schotzau1999hp, schotzau2000hp} to solve parabolic problems and extended in \cite{feistauer2007space} to a linear nonstationary convection-diffusion-reaction problem. This method, denoted CG($p$)DG($q$), used a piecewise polynomials of degree $p$ and $q$ respectively for the space and time discretization.  Feistauer et al. \cite{feistauer2011analysis} proposed the theory of error estimates for  CG($p$)DG($q$) applied to a nonstationary convection-diffusion problem with a nonlinear convection and linear diffusion. The DG method is used in both space and time  by Feistauer et al. \cite{vcesenek2012theory, feistauer2011, dolejvsi2015discontinuous} to solve the nonstationary parabolic problems with nonlinear convection and diffusion.

Other than the DG time discretization, L{\"o}rcher et al.  \cite{lorcher2007dis} used the LTS method (denoted ADER-DG) based on arbitrary high-order derivatives methods and allowed every element of the mesh to have its own time-step, which is dictated by the element size. The  ADER-DG schemes, as presented for electromagnetism \cite{taube2009high} and elastic wave propagation in \cite{dumbser2007arbi}, were obtained by the extension to the DG framework, of the ADER finite volume (ADER-FV) approach which was developed by Toro et al. \cite{titarev2002ader}. The ADER-DG scheme was used by Fambri et al. \cite{fambri2017space} on space-time adaptive meshes for compressible Navier-Stokes equations and the equations of viscous and resistive magnetohydrodynamics in two and three space-dimensions. 
Angulo et al. \cite{angulo2014causal} introduced the LTS schemes (denoted LTS-LF) based on the leapfrog (LF)  and Runge-Kutta (RK) time integrators, where the mesh is sorted  into different sub-domains with appropriate time step on each. In 2009, Diaz and Grote introduced an energy-conserving LTS-LF scheme \cite{diaz2009energy} for the acoustic wave equation, which they extended in 2015 into a multi-level version \cite{diaz2015multi}.  Rietmann et al. \cite{rietmann2015load} developed a new LTS method based on the Newmark scheme for large scale wave propagation, which also can be extended to accommodate multiple sub-domains of mesh refinement.

Unfortunately, the DG method in space and time considered in  \cite{vcesenek2012theory, feistauer2011, dolejvsi2015discontinuous, lorcher2007dis, angulo2014causal, diaz2009energy, diaz2015multi, rietmann2015load} were still special, since they always had boundaries in time aligned with the time direction i.e. the spatial boundaries are independent of the time. To overcome this limitation, an alternative space-time DG method was introduced by  van der Vegt et al. in \cite{van11, van22} for inviscid compressible flows and extended in \cite{van44} to the compressible Navier-Stokes equation. The key feature of this space-time DG method is that no distinction is made between space and time variables and the DG discretization is directly and simultaneously performed in space and time. This then provides  flexibility to deal with time dependent boundaries, deforming elements and naturally results in a conservative discretization, even on deforming, locally refined meshes with hanging nodes. A complete hp-error and stability analysis of the space-time DG discretization for the linear advection-diffusion equation is given in \cite{van33}.

However, the bottleneck of these LTS methods, based on the DG discretization, is that they lead to a large discrete problem in space-time, especially in the presence of complex geometry or localized small-scale physics. As a consequence, they can become very expensive in terms of storage and computational time. Thus, we focus on splitting the solution domain on several small regions, yielding several low dimension system of ODEs from the DG spatial discretization, which can then be solved separately. We look at two
approaches: the first based on \cite{LTSP4, Wang6, VAB17, dryja1991, laevsky1993, vab94parallel} where the sub-domains are overlapped and the second using non-overlapping sub-domains based on \cite{LTSP9, dawson1991finite, dawson1992explicit}. Let us denote by GTS-DG, the global time stepping solver of the DAREs that combines the DG method and a time integrator with unique time step on the entire solution domain. 

The format of this paper is as follows.  In  \sectref{\ref{LTS-DG-sc}},  we describe our proposed  LTS-DG schemes. Various numerical results are presented in \sectref{\ref{ALL_num_res}} to validate and compare the efficiency of the proposed LTS-DG schemes against GTS-DG schemes, and the conclusion follows in \sectref{\ref{sec:conclusions}}.

\section{LTS-DG schemes\label{LTS-DG-sc}}
\begin{figure}
\centering
\subfloat[]{\includegraphics[width=0.45\textwidth,height=0.15\textheight]{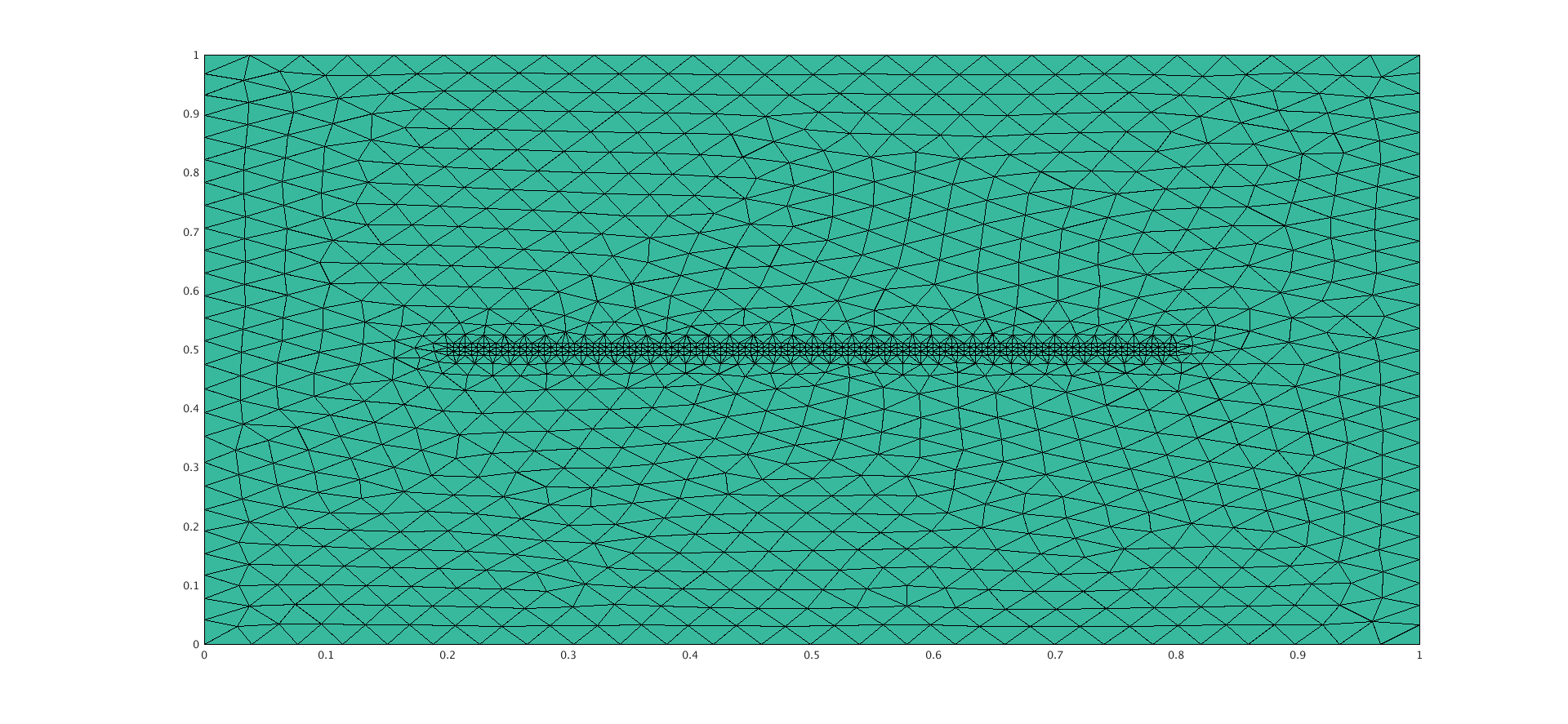}}
\subfloat[]{\includegraphics[width=0.45\textwidth,height=0.15\textheight]{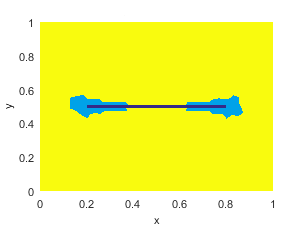}}
\caption{(a) Unstructured mesh of domain with fracture using $distmesh2d$.
(b) Sub-domains: $\Omega_0$ (black), $\Omega_1$ (blue) and $\Omega_2$ (yellow).
}
\label{frac}
\end{figure}

We construct our LTS-DG schemes by combining domain decomposition
techniques, the DG spatial discretization and the time
integrators. The approach follows three basic steps: 
\begin{enumerate}
\item We use a priori knowledge of the local behaviour of the solution
  due to   fractures, walls, corners, obstacles, point loads, etc, to
  construct a refined mesh of the spatial domain. For example, this is
  illustrated in \figref{\ref{frac}} (a) by showing the refined mesh for the domain with fracture.

\item We choose the local time step on each element of the mesh such
  that it is proportional to the element size (similar to
  \cite{lorcher2007dis}) and inversely proportional to the norm of the
  fluid velocity $\textbf{v}$ (if it is different from zero). This
  splits the solution domain, $\Omega$, into sub-domains, $\Omega_i$
  with the local time steps $\Delta t_i$ for all
  $i\in\{0,\cdots,m\}$. 
\item Finally, we use either interpolation or extrapolation techniques to estimate the solution at the internal boundary, $\Gamma_i$, given by 
\begin{equation}
\Gamma_i = \partial\Omega_i \setminus \partial\Omega,
\end{equation}
for all $i\in\{0,\cdots,m\}$. So, one can solve the PDEs independently on each sub-domains, $\Omega_i$, using the DG method combined with the standard time integrators and local time step $\Delta t_i$. 
\end{enumerate}

We now look at these steps in more details. The first step can be implemented, for example, using the MATLAB's code $distmesh$ \cite{Persson04, Persson05}. We complete the second step by assuming that, beside being proportional to the element size and inversely proportional to the norm of the fluid velocity, the local time step on the sub-domain $\Omega_i$ is given by 
\begin{equation}\label{lamd_dt}
\Delta t_i = \lambda_i\times\Delta t,
\end{equation}
with $\lambda_i\in\mathbb{N}$, for a given $\Delta t\in \mathbb{R}^+$
and all $i\in\{0,\cdots,m\}$. Thus, for a given initial time $T^0$ and
final time $T^1$, we have a time synchronization over all sub-domains
at some time $t\in [T^0, T^1]$. 
We illustrate how this can be achieved, in two dimension case,
ensuring that the CFL condition is satisfied on each element $E$ of
the triangulation. Let
\begin{equation}\label{local_time1}
\Delta t_E = 2^{-\mathbf{N}^E}, \;\;\;\;\mathbf{N}^E=\lceil\log_2\left(  \frac{\parallel \textbf{v}^E \parallel }{R_E\mathbf{C}_{max}} \right) \rceil,
\end{equation}
where $\lceil\cdot\rceil$ is the ceiling function \cite{ceil10}, $\mathbf{C}_{max}\in \mathbb{R}$ is the Courant number \cite{CFL1}, $R_E$ is the radius of the incircle of the element $E$. This ensures that the CFL condition is verified on each element $E$. \equaref{\ref{local_time1}} yields the decomposition of the global solution domain, $\Omega$, into sub-domains $\Omega_i, i\in\{0,\cdots,m\}$ defined as follows 
\begin{equation}
\Omega_i = \lbrace E\in \mathcal{T}_h \text{ such that } \Delta t_E = \Delta t_i\rbrace,
\end{equation} 
where $\Delta t_i$ is a fixed value in $\mathbb{R}$ for a given value of $i$ in $\{0,\cdots,m\}$. The synchronized time $t = t^n$ while advancing locally the solution from $T^0$ are given by 
\begin{equation}
t^n = T^0 + n\times \Delta t_{\max}, \; \text{ with }\;\Delta t_{\max} = \max\{ \Delta t_i, i=0,\cdots,m \}.
\end{equation}

%
The sub-domains, obtained when we applied the second step to the fracture problem, are illustrated in \figref{\ref{frac}}(b). Note from \figref{\ref{frac}}(b) that $m=2$ with the sub-domains $\Omega_0, \Omega_1$ and $\Omega_2$ respectively represented by the color black, blue and yellow. One can also see that the sub-domain $\Omega_1$ is the union of two disjoint regions while its interior boundary $\Gamma_1$ is shared with the sub-domains $\Omega_0$ and $\Omega_2$.
%
Once the sub-domains are defined, the large time-dependant system of ODEs, obtained from the DG space discretization of the DAREs on the $\Omega$, can then be split into $m+1$ smaller systems of ODEs, denoted $\text{SO}_i$ for all $i\in\{0,\cdots,m\}$,  given by
\begin{equation}\label{equa_locl}
(\text{SO}_i):\;\;\mathbb{M}_i\frac{d}{dt}X_i + \mathbb{S}_iX_i = \mathbb{F}_i  + \mathbb{B}_i + \mathbb{S}_i^eX_i\Bigr|_{\Gamma_i}\;\text{ on } \;\Omega_i\times[0,T].
\end{equation}
Here, $X_i$, $\mathbb{S}_i, \mathbb{M}_i$, $\mathbb{F}_i$ and $\mathbb{B}_i$  respectively represent the local solution, stiffness matrix, the mass matrix, source term and the contribution of the global boundary condition on the sub-domain $\Omega_i$. The matrix $\mathbb{S}_i^e$ is used to weakly enforce the internal boundary condition. Therefore, the sum of the last two terms at the right hand side of \equaref{\ref{equa_locl}} enforces the local boundary condition on $\partial\Omega_i$. Because the interior penalty discontinuous galerkin (IP-DG) method is compact  and the global mass matrix obtained from the IP-DG spatial discretization is either block-diagonal or diagonal, see \cite{SP2}, then the local entities $X_i$, $\mathbb{S}_i, \mathbb{M}_i$, $\mathbb{S}_i^e$, $\mathbb{F}_i$ and $\mathbb{B}_i$ can be easily extracted from their global values. An example of this extraction is shown in \subsectref{\ref{OETO_expe}}.

Let us introduce the following notation
\begin{equation}
t_{i}^j = T^0 +j\times \Delta t_i,\qquad s_i = \frac{T^1 - T^0}{\Delta t_i},\qquad
X_i^j = X_i\Bigr|_{t=t_{i}^j}, 
\end{equation}
for all $i\in\{0,\cdots,m\}$ and all $j\in\{0, \cdots, s_i\}$. If one can estimate $X_i^j\Bigr|_{\Gamma_i}$, by any means, then the local solution $X_i^j$ at time $t_{i}^j$ can be obtained from its initial value $X_i^0$, by applying the time integrator schemes to the local system $\text{SO}_i$ given by \Equaref{\ref{equa_locl}}, with the uniform local time step $\Delta t_i$.

As a consequence, the construction of our LTS-DG schemes is reduced to finding a way to estimate the component $X_i^j\Bigr|_{\Gamma_i}$ at any time $t_i^j$ for all $i\in\{0,\cdots,m\}$ and all $j\in\{0, \cdots, s_i\}$. We only need to describe the LTS-DG schemes to advance the global solution $X$, on the solution domain $\Omega$, from its known value at the synchronized time $t^n$ (as it happens at the initial time $T^0$) to the synchronized time $t^{n+1}$. The process can then be repeated, in order to estimate the global solution $X$ at the final time $T^1$ from the global solution $X^0$ at the initial time $T^0$. Next, in  \subsectref{\ref{Over_ref}} and \subsectref{\ref{NOver_ref}}, we described two different techniques respectively refer as overlap LTS-DG and non-overlap LTS-DG methods to locally advance the solution from synchronized time $t^n$ to $t^{n+1}$. 

\subsection{Overlap LTS-DG schemes (OLTS-DG)\label{Over_ref}}
The key idea of the proposed overlap LTS-DG methods, denoted OLTS-DG, is to overlap the different sub-domains obtained from the decomposition of the solution domain $\Omega$, in order to extrapolate the components $X_i^j\Bigr|_{\Gamma_i}$ at any time $t_i^j$ for all $i\in\{0,\cdots,m\}$ and all $j\in\{0, \cdots, s_i\}$. 
This approach appeared in \cite{LTSP4}, where a Crank-Nicolson scheme was used for the time-space discretization of the one spatial dimension heat equation. The authors proved that without local refinement in time ($\Delta t_i = \Delta t$) and space ($h_i = h$) this scheme is stable, provided that
\begin{equation*}
\Delta t \leq C\left( \frac{L}{\log L}  \right)^2h^2,
\end{equation*}
where $Lh$, for $L\in \mathbb{N}$, is the size of the overlap, and an error estimate of the form $O(\Delta t^2 +h^2)$. So, in this case, increasing the size of the overlap can reduce the stability constraint
on the time step. To avoid the stability constraint, Ewing et al. \cite{ewing1994finite} used a standard centred finite difference scheme in space with backward Euler in time for a linear DAREs. More recently, an approach based on domain decomposition and finite volume discretization, has been proposed by Faille et al. \cite{faille2009two} for the one dimensional heat equation.  It was used by Gander et al. \cite{gander2005overlapping} to investigate the one dimensional convection dominated nonlinear conservation laws. 

Here, we introduce new schemes that extend this approach to the DAREs
in one, two or three spatial dimensions. We use the DG method for the
space discretization and time integrators such as Impl, ETD, EXPR for
the resolution in time of \equaref{\ref{equa_locl}}, in order to avoid
numerical instability.
\subsubsection{Overlapping procedure of the domain solution}
Once the sub-domains $\Omega_i,\;i\in\{0,\cdots,m \}$ are obtained, we
overlap them by pushing the internal boundary $\Gamma_i$ in the
direction of the outward normal vector. 
\newpage
\begin{figure}
\centering
\includegraphics[height=0.3\textheight,width=0.9\textwidth]{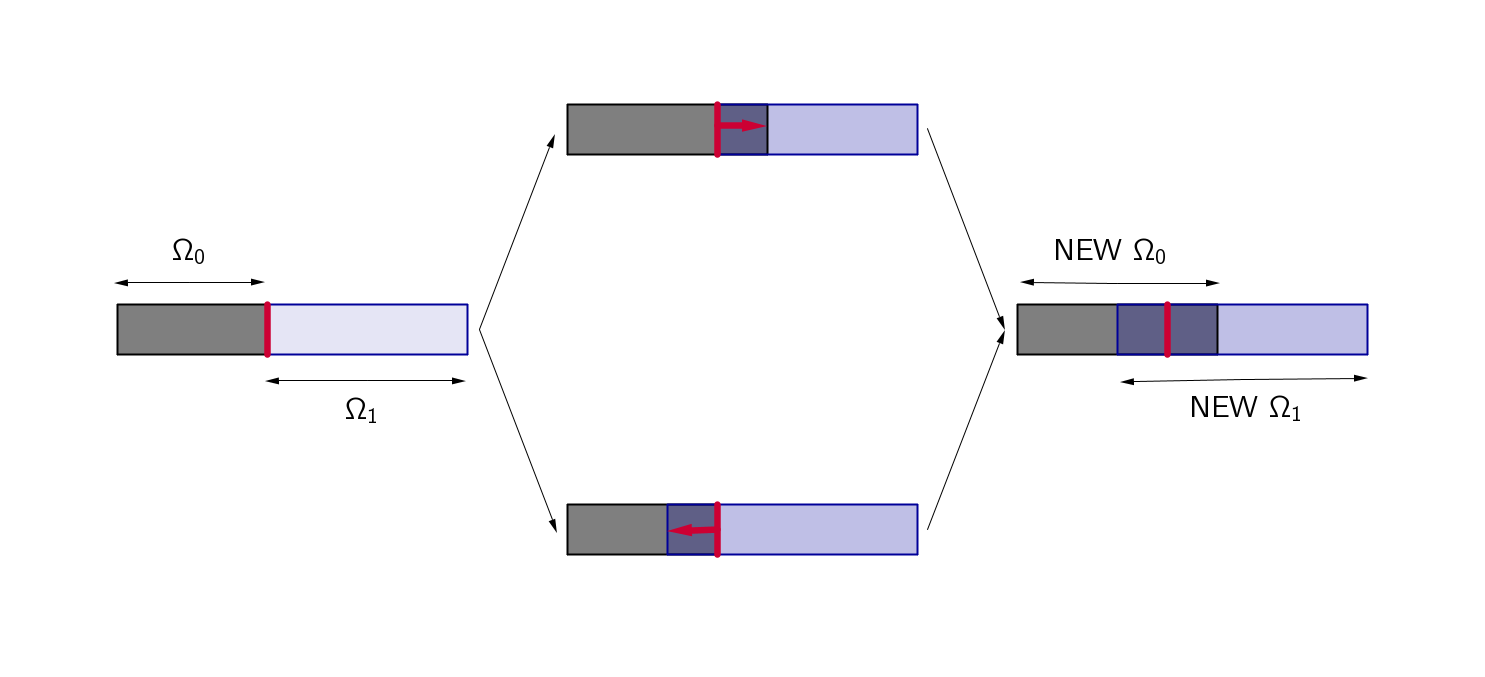}
\caption{The procedure to overlap the regions for $\Omega = \Omega_0\cup \Omega_1$. On the left: original sub-domains $\Omega_0$ and $\Omega_1$. Middle: we push the internal boundary in the direction of the outward nornal vector. Right: new sub-domains $\Omega_0$ and $\Omega_1$\label{overlap_metho}}
\end{figure}
During the overlapping procedure, if any new sub-domain $\Omega_i$ swallows entirely another initial sub-domain $\Omega_j$, then we set $\Delta t_j=\Delta t_i$ so that the sub-domain $\Omega_j$ will be included in $\Omega_i$. Later on, in \subsectref{\ref{OGAT_expe}}, we investigate numerically the effect of the size of the overlap on the accuracy of our OLTS-DG schemes.

Let us denote $\Gamma_{i,j}$ the part of the internal boundary
$\Gamma_i$ included in the sub-domain $\Omega_j$ with $i\neq j$
(i.e. $\Gamma_{i,j} = \Gamma_{i}\cap \Omega_j)$. 
Every time we advance the local  solution
$X_j$  on  $\Omega_j$, we can update the $X_i\Bigr|_{\Gamma_{i,j}}$.

For a given $r\in \{1, \cdots,\frac{\Delta t_{\max}}{\Delta t}\}$, the set of eligible sub-domain $S^r$ is the set of sub-domains on which the known solution has to be advanced locally to the time $t^{n_r} = t^n + r\times \Delta t$. Then we have
\begin{equation}
S^r = \left\lbrace \Omega_i \Bigr|\exists j\in \mathbb{N}, t_i^j = t^{n_r}    \right\rbrace.
\end{equation}

Note that there is a freedom in the order of which the sub-domains are updated. For example, it could either be in the increasing or decreasing order of local time step. If the time integrator $\text{INT}\in\{$Impl, ETD, EXPR $\}$ is used to advance locally the solution, we denote $DG_{\text{OLTSD-INT}}$ and $DG_{\text{OLTSI-INT}}$ the OLTS-DG scheme that updates the solution on the eligible sub-domains $S^{r}$ respectively in the decreasing and increasing order of the local time step. In \subsectref{\ref{OGAT_expe}}, we compare the accuracy of the  $DG_{\text{OLTSD-Impl}}$ and $DG_{\text{OLTSI-Impl}}$ and examine how the direction of the bulk velocity of the DAREs or the size of the overlap affect their accuracy. Unless stated, the OLTS-DG method considered for the numerical experiments is the $DG_{\text{OLTSD-INT}}$. Next, we describe step by step the algorithm of the OLTS-DG scheme, $DG_{\text{OLTSD-INT}}$, in order to advance the solution from a synchronized time $t^n$ to $t^{n+1}$.


\subsubsection{Description of the OLTS-DG algorithm}
\begin{figure}[H]
\centering
\includegraphics[height=0.2\textheight,width=0.9\textwidth]{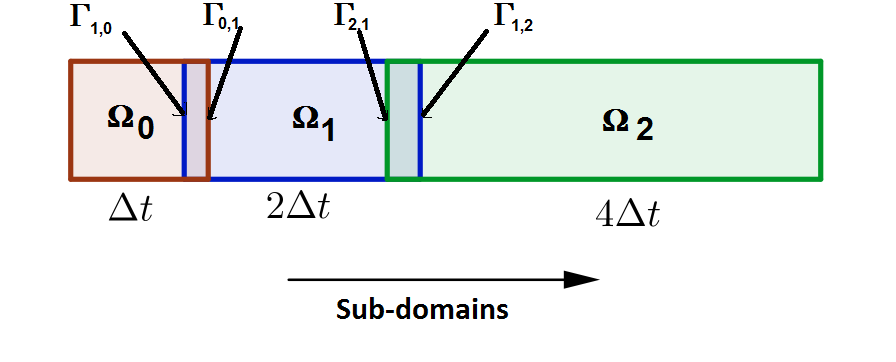}
\caption{Overlapped sub-domains of the solution domain.\label{overlap_metho1}}
\end{figure}

In this section, we describe step by step the algorithm of the OLTS-DG scheme, $DG_{\text{OLTSD-INT}}$, in order to advance the solution from a synchronized time $t^n$ to $t^{n+1}$. To that end, we consider the overlapped sub-domains illustrated in \figref{\ref{overlap_metho1}}, where the solution domain is split into three different sub-domains (i.e.  $m=2$) with the coefficient $\lambda_i$ in \Equaref{\ref{lamd_dt}} given by $\lambda_i = 2^i$ for all $i=0,\cdots,m$.
The overlapped sub-domains $\Omega_0$, $\Omega_1$ and $\Omega_2$ are respectively represented by the color red, blue and green. Note from \figref{\ref{overlap_metho1}} that the internal boundaries are given by $\Gamma_0 = \Gamma_{0,1}, \; \Gamma_1 = \Gamma_{1,0}\sqcup \Gamma_{1,2},\;\; \Gamma_2 = \Gamma_{2,1}$.


\newpage
For a given time $t^{n_r}$, let us denote $\mathbb{X}_i^{n_r}$ and $\mathbb{X}_{i,j}^{n_r}$ the restriction of the global solution $X$  respectively to the internal boundaries $\Gamma_i$ and $\Gamma_{i,j}$ i.e.
\begin{equation}
\mathbb{X}_i^{n_r} = X\Bigr|_{\Gamma_i,t = t^{n_r}},\;\;\mathbb{X}_{i,j}^{n_r} = X\Bigr|_{\Gamma_{i,j},t = t^{n_r}}.
\end{equation}
Now let discuss step by step, how to update the solution on these interior boundaries in order to advance the solution from the synchronized time $t_n$ to $t_{n+1}$.

\begin{itemize}
\item \textbf{Step one:} for $r=1$, the set of eligible sub-domains is
  given by $S^1 = \{{\Omega_0}\}$  and requires $\mathbb{X}_0^{n_1}$
  to advance locally on $\Omega_0$ from $t^n$ to $t^{n_1} = t^n +
  \Delta t$. To that end, we then use the extrapolation
  $\mathbb{X}_0^{n_1} = \mathbb{X}_0^{n}$. The completion of this step
  defines $\mathbb{X}_{1,0}^{n_1}$, as illustrated in
  \figref{\ref{overlap_metho23}} (a).
\begin{figure}[H]
\centering
\subfloat[]{\includegraphics[height=0.35\textheight,width=0.48\textwidth]{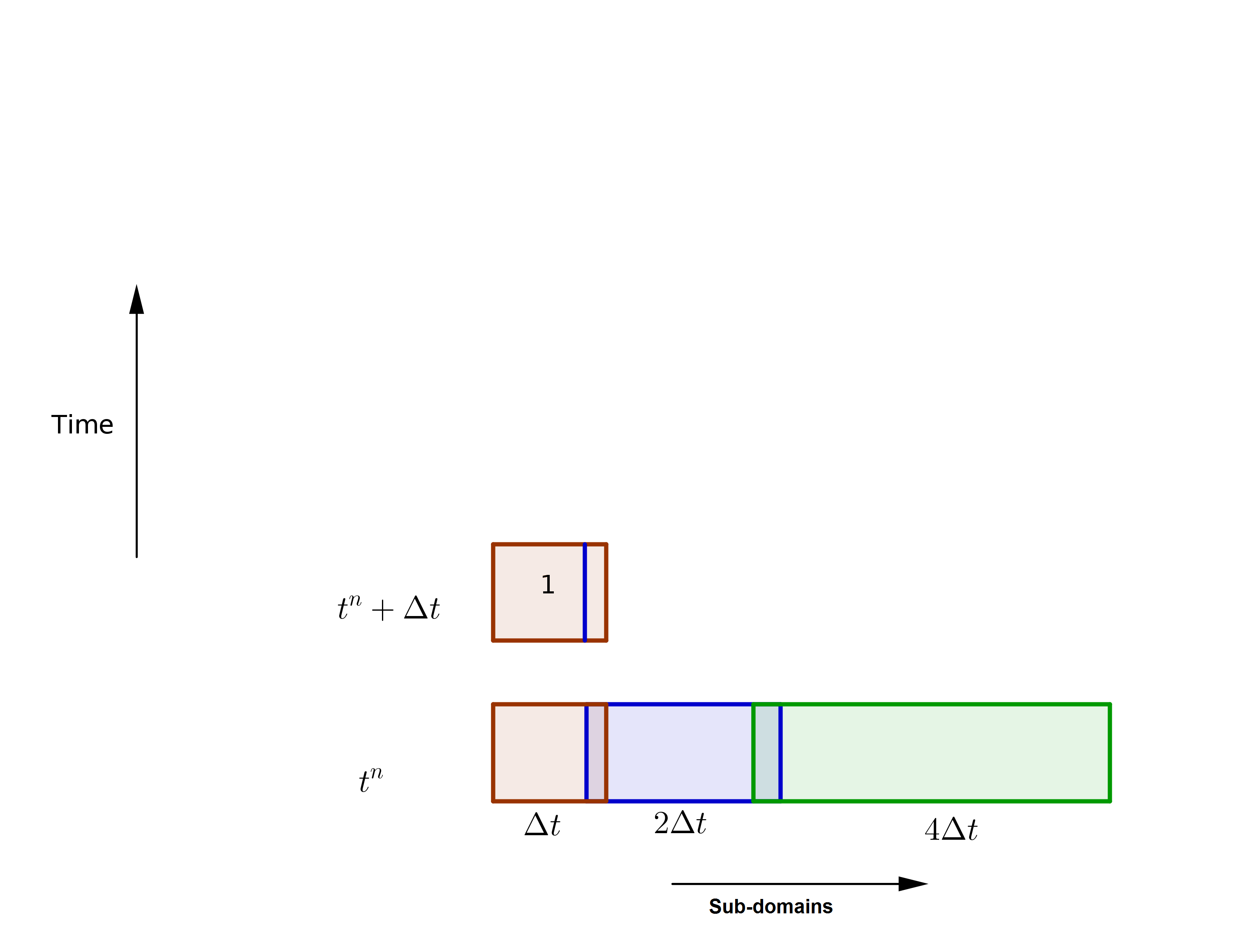}}
\subfloat[]{\includegraphics[height=0.35\textheight,width=0.48\textwidth]{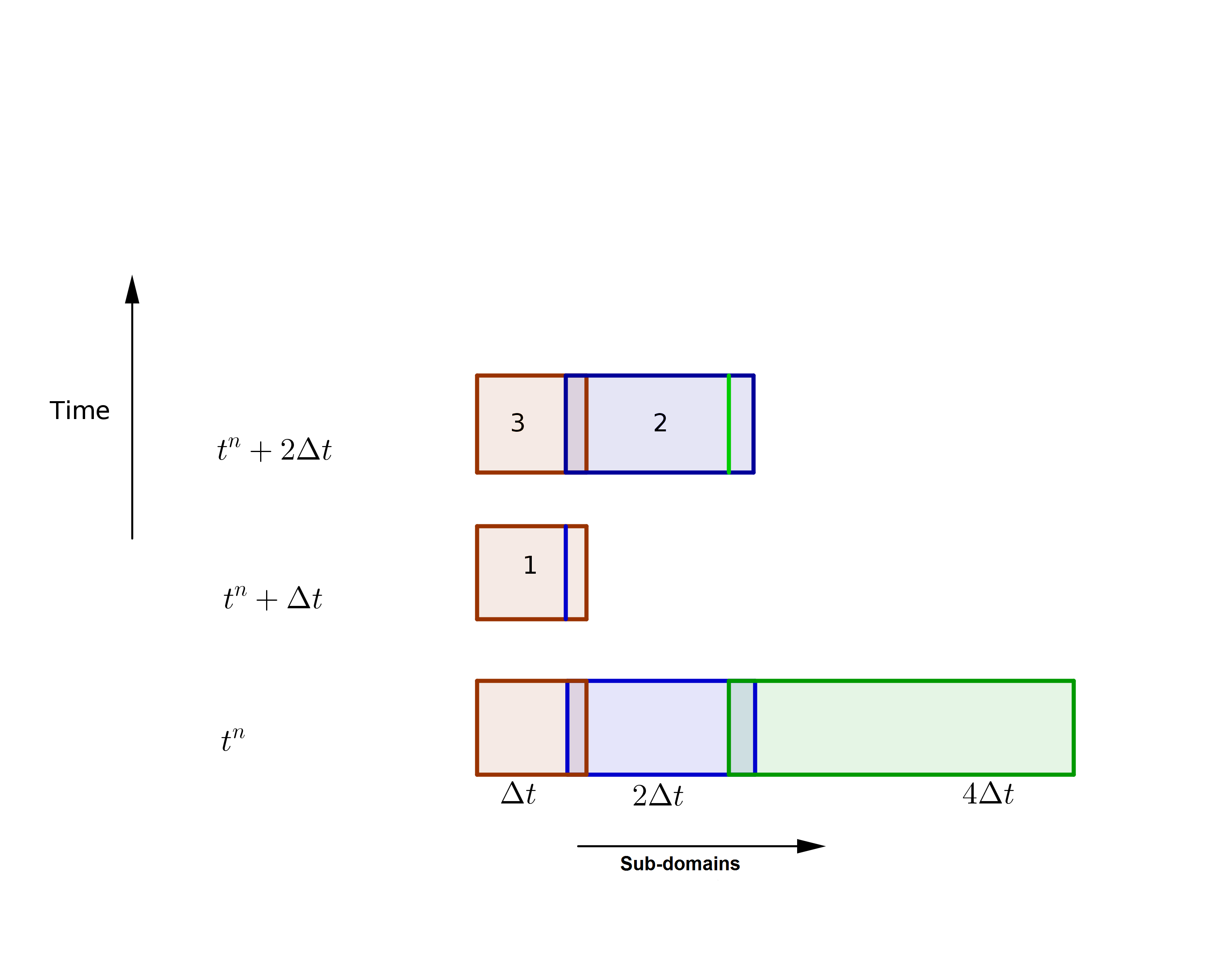}}
\caption{(a) The eligible solution advanced to
  $t^{n_1}$ and (b) to $t^{n_2}$.}
\label{overlap_metho23}

\end{figure}

\item \textbf{Step two:} for $r=2$, the set of eligible sub-domains is given by  $S^2 = \{\Omega_0,\Omega_1\}$, and we require $\mathbb{X}_0^{n_2}$ and $\mathbb{X}_1^{n_2}$ to locally advance the solution to $t^{n_2}= t^n + 2\times\Delta t$. So, we first use the extrapolations $\mathbb{X}_{1,0}^{n_2} = \mathbb{X}_{1,0}^{n_1}$ and $\mathbb{X}_{1,2}^{n_2} = \mathbb{X}_{1,2}^{n}$ to locally advance the solution on $\Omega_1$ from $t^n$ to $t^{n_2}$. Note from \figref{\ref{overlap_metho23}}(b) that the completion of this simulation defines $\mathbb{X}_0^{n_2}$ and $\mathbb{X}_2^{n_2}$. We finally use $\mathbb{X}_0^{n_2}$ to advance locally the solution on $\Omega_0$ from $t^{n_1}$ to $t = t^{n_2}$. 

\item \textbf{Step three:} for $r=3$, the set of eligible sub-domains
  is given by $S^3 = \{{\Omega_0}\}$  and require $\mathbb{X}_0^{n_3}$
  to advance locally on $\Omega_0$  to $t^{n_{3}}= t^n + 3\times\Delta
  t$. We use the extrapolation $\mathbb{X}_0^{n_3} =
  \mathbb{X}_0^{n_2}$. This is illustrated in
  \figref{\ref{overlap_metho45}}(a)  and it shows that the completion of this
  step defines $\mathbb{X}_{1,0}^{n_3}$.
\begin{figure}[H]

\centering
\subfloat[]{\includegraphics[height=0.35\textheight,width=0.48\textwidth]{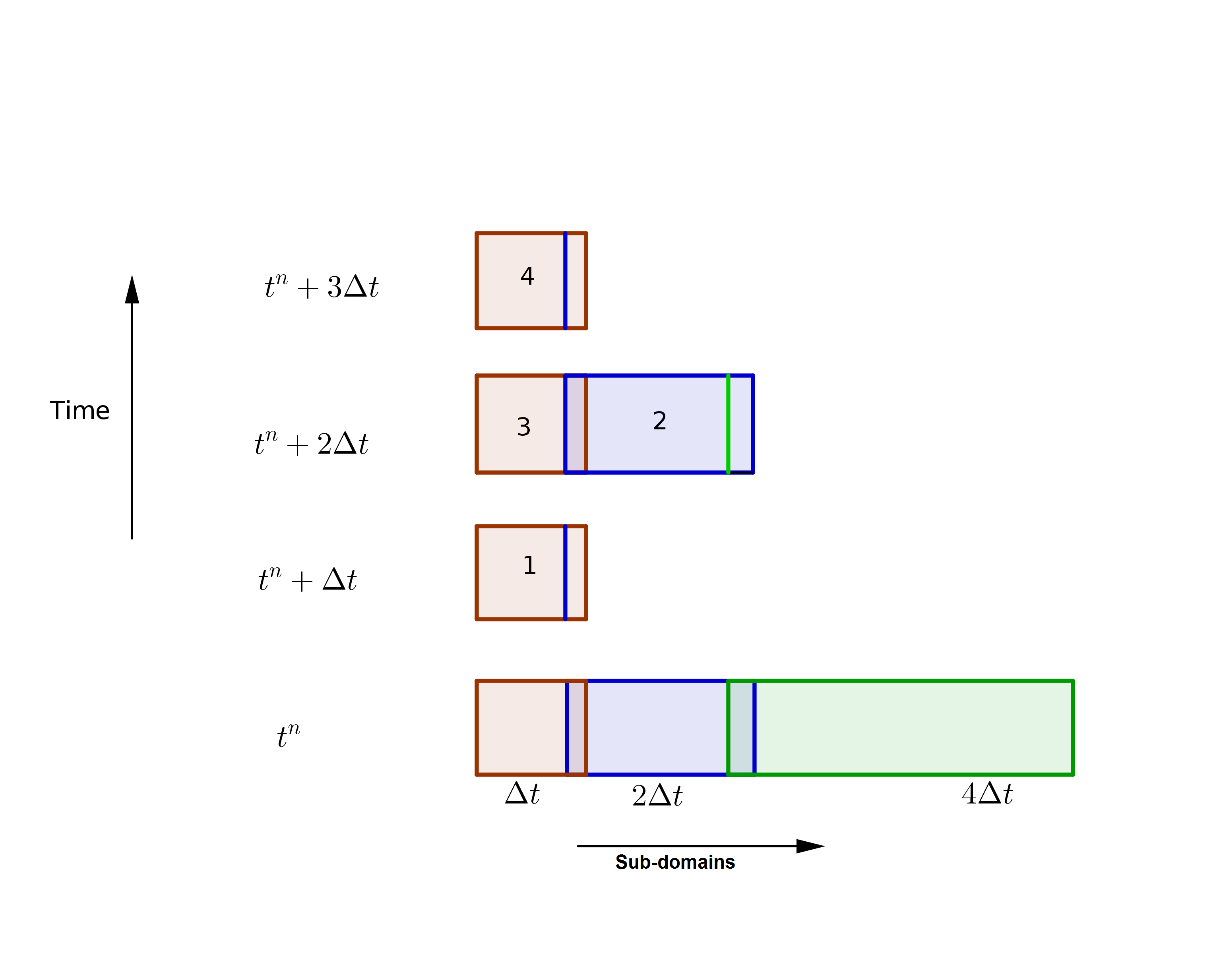}}
\subfloat[]{\includegraphics[height=0.35\textheight,width=0.48\textwidth]{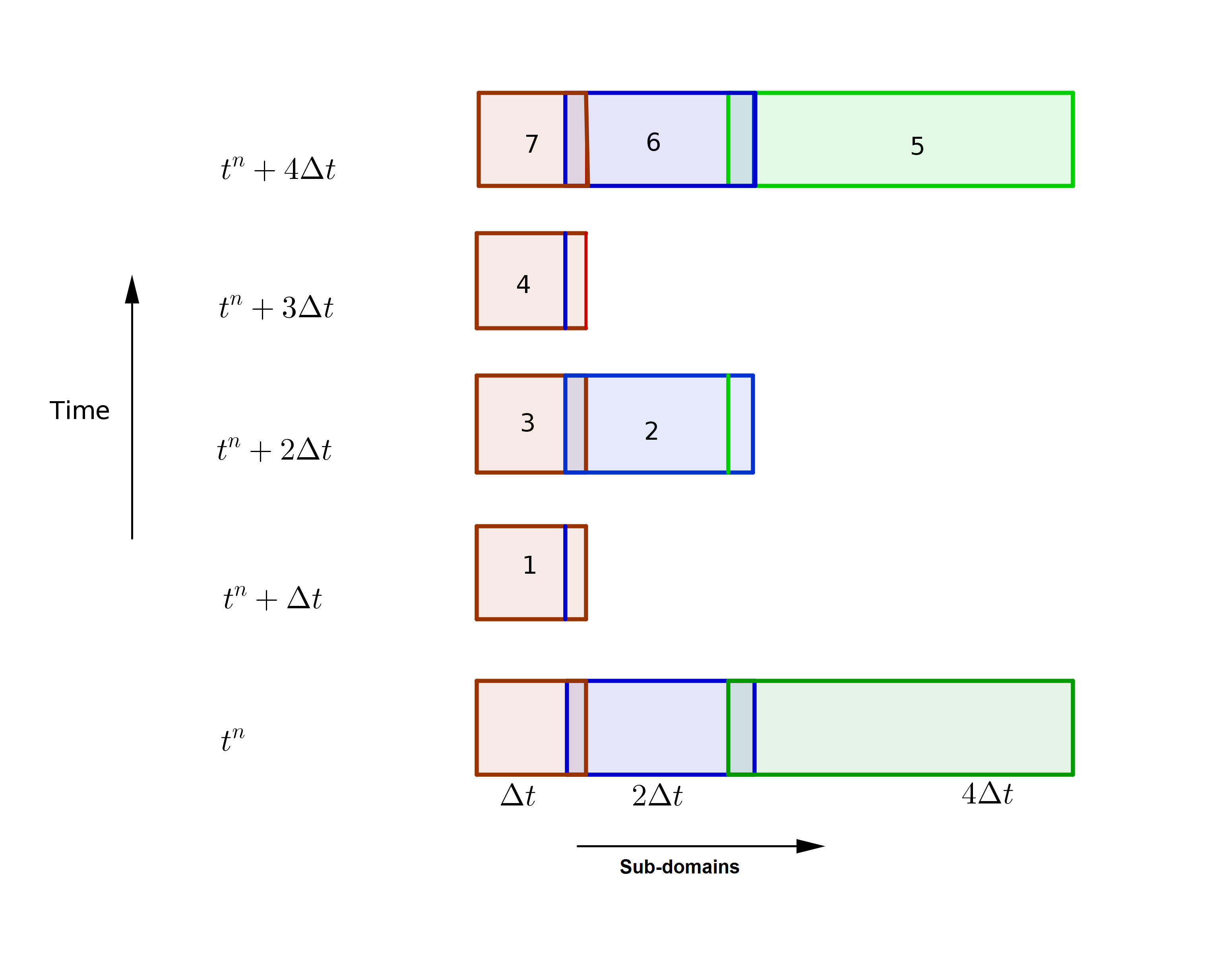}}
\caption{The eligible solution advanced to $t^{n_3}$ in (a)
  and $t^{n_4}$ in (b).}
\label{overlap_metho45}
\end{figure} 

\item \textbf{Step four:} for $r=4$, the set of eligible sub-domains
  is given by  $S^4 = \{\Omega_0,\Omega_1, \Omega_2\}$, and we require
  $\mathbb{X}_0^{n_4}, \mathbb{X}_1^{n_4}$ and  $\mathbb{X}_2^{n_4}$
  to locally advance the solution to $t^{n_4}= t^n + 4\times\Delta
  t$. We first locally advance the solution on $\Omega_2$ from $t^n$
  to $t= t^{n_4}$, using the extrapolation $\mathbb{X}_2^{n_4} =
  \mathbb{X}_2^{n_2}$. This is illustrated in \figref{\ref{overlap_metho45}}
  (b) and it shows that this simulation defines
  $\mathbb{X}_{1,2}^{n_4}$. We then use $\mathbb{X}_{1,2}^{n_4}$
  obtained and the extrapolation $\mathbb{X}_{1,0}^{n_4} =
  \mathbb{X}_{1,0}^{n_3}$, to locally advance the solution on
  $\Omega_1$ from $t^{n_2}$ to $t^{n_4}$. This is also illustrated in
  \figref{\ref{overlap_metho45}} (b). It shows that this simulation  defines
  $\mathbb{X}_0^{n_4}$, which we then use to locally advance on
  $\Omega_0$ from $t^{n_3}$ to $t^{n_4}$. 
\end{itemize}

At this point, the time is synchronized across the whole domain $\Omega$. By repeating this process (i.e. from step one to step four) $m^s = \frac{T^1 - T^0}{\Delta t_{\max}}$ times, we can estimate the solution at the final time $T^1$, from the solution at the initial time $T^0$. One can implement \Algref{\ref{Over_pseudo_algo}}, where  
 $S$ is the set of all overlapped sub-domains $\Omega_i$,
$\mathbf{E}$ and $\mathbf{A}_i$ respectively represent the extrapolation procedure and the iterative function of the standard time integrators used to solve the local system $\text{SO}_i$.

%
%
%
%
%
%
%
%
%
%
%
%
%
%
%
%
\newpage
\begin{algorithm}[H]
\caption{Pseudo algorithm of the overlap LTS method. \label{Over_pseudo_algo}}
\begin{algorithmic}
\For{$i= 0 \text { to } m-1$}\Comment{initialize the time of the known local solution $X_i$}
    \State $t^{known}_i=T^0$
\EndFor
\For{$n= 0 \text { to }  m^s-1$}\Comment{Loop to advance the solution from $T^0$ to $T^{1}$ }
    \State $t^n = T^0+j\times \Delta t_{\max}$ \Comment{Compute the synchronized time $t^{n}$}
   \For{$r = 1 \text { to }  \frac{\Delta t_{\max}}{\Delta t}$}\Comment{Loop to advance the solution from $t^n$ to $t^{n+1}$}
    	\State $t^{n_r} = t^n + r\times \Delta t$ \Comment{Compute the local advancing time $t^{n_r}$}
    	\For{$i= 0 \text { to }  m-1$}\Comment{Loop to advance locally the solution on the eligible sub-domains to $t^{n_r}$}
    		\If{$\mod(r\times \Delta t, \Delta t_i)=0$}\Comment{If $\Omega_i$ is eligible at the time $t^{n_r}$, then advance locally}
    			\State $\mathbb{X}_i^{n_r} = \mathbf{E}\left(  X_{l_i} \Bigr|_{ t^{known}_{l_i}}, \; \forall l\neq i\rbrace\right) $ \Comment{Extrapolate the local solution $X_i$ at the boundary $\Gamma_i$ to $t^{n_r}$}
    			\State $ X_i \Bigr|_{t^{n_r}} = \mathbf{A}_i\left(X_i \Bigr|_{t^{known}_i},  \mathbb{X}_i^{n_r}  \right) $\Comment{Advance locally the solution on $\Omega_i$ to $t^{n_r}$}
				\State $t^{known}_i = t^{n_r}$\Comment{Update the local time $t^{known}_i$ of the known local solution $X_i$}
			\EndIf
		\EndFor
	\EndFor
\EndFor
\end{algorithmic}
\end{algorithm}

\subsection{Non overlap LTS-DG schemes (NOLTS-DG)\label{NOver_ref}}

Another way to explicitly estimate the value of $X\Bigr|_{\Gamma_i}$ appeared in \cite{LTSP9} in finite difference context for heat equation. In which case there is no need of extending the boundary of the sub-domain $\Omega_i$, once the local time steps $\Delta t_i$ are defined. The key idea of the non overlap method, NOLTS-DG, is to first advance the solution globally to the time $t^n+\Delta t^*$ from the known solution at time $t^n$, where the global time step $\Delta t^*$ larger than the maximum local time step $\Delta t_{max}$. This step is called the prediction step and is followed by an interpolation  to obtain the values of $X\Bigr|_{\Gamma_i}$ needed to advance the solution locally on the sub-domain $\Omega_i$. This last step is called the correction step. It has been applied in finite element context \cite{dawson1991finite} and discontinuous Galerkin context \cite{dawson1992explicit} for parabolic equations.

We now develop new schemes that extend this approach to the DAREs in one, two or three spatial dimensions, using the DG method for the space discretization and time integrators such as Impl, ETD or EXPR for the resolution of the local system $\text{SO}_i$. 

\subsubsection{Non overlap LTS-DG algorithm}
Once again, we consider the case where the solution domain $\Omega$ is split into three different sub-domains $\Omega_i$ with the local time step $\Delta t_i = 2^{i}\times\Delta t$ for all $i=0,\cdots,m$, $m=2$ and a given time step $\Delta t$. Therefore, in order to obtain the component $X$ of the concentration entirely on $\Omega$ at the time $t^{n+1}$, from its known value at the time $t^n$ using the non overlap LTS method, we use the following steps.

\begin{figure}[H]
\centering
\subfloat[]{\includegraphics[height=0.3\textheight,width=0.48\textwidth]{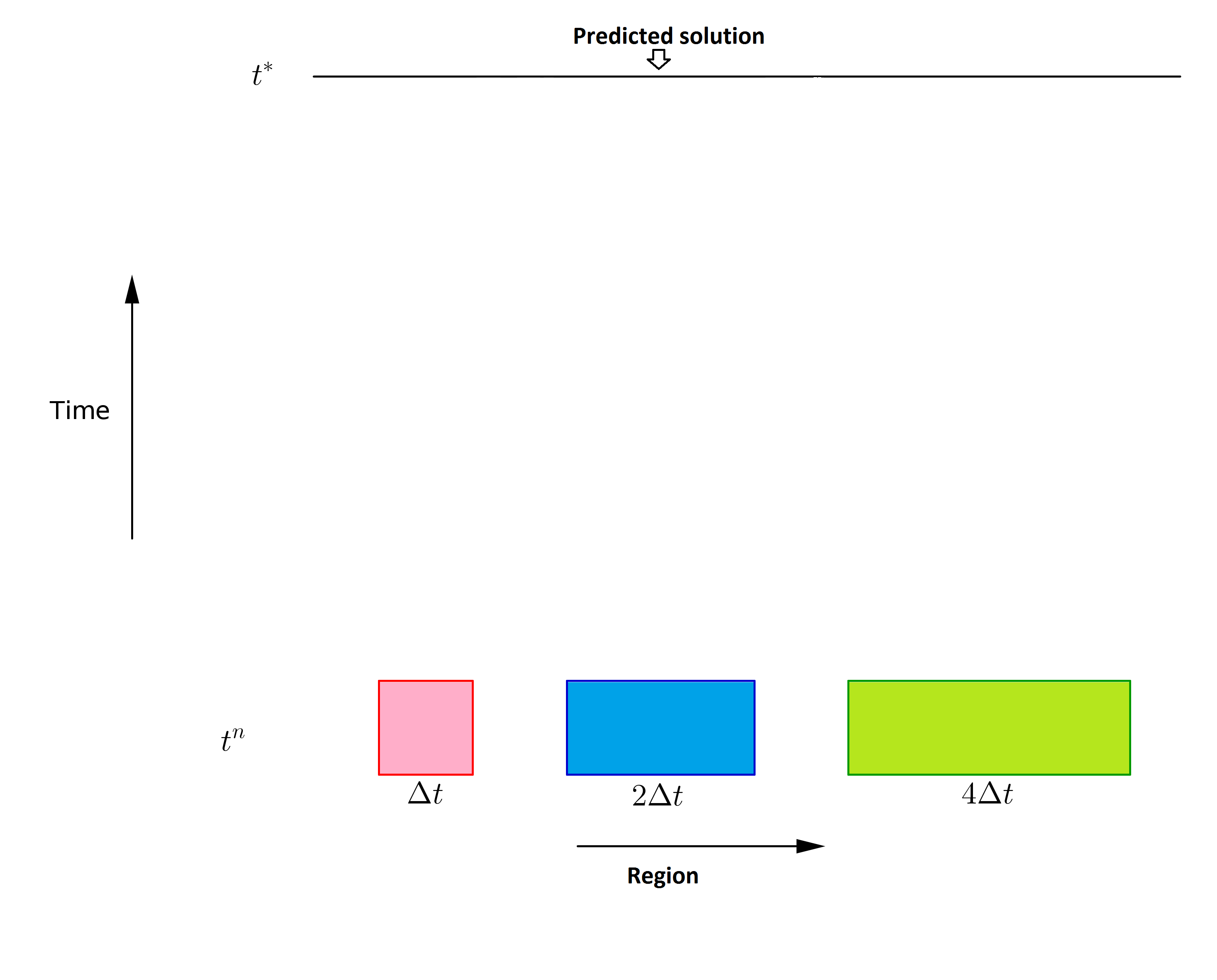}}
\subfloat[]{\includegraphics[height=0.3\textheight,width=0.48\textwidth]{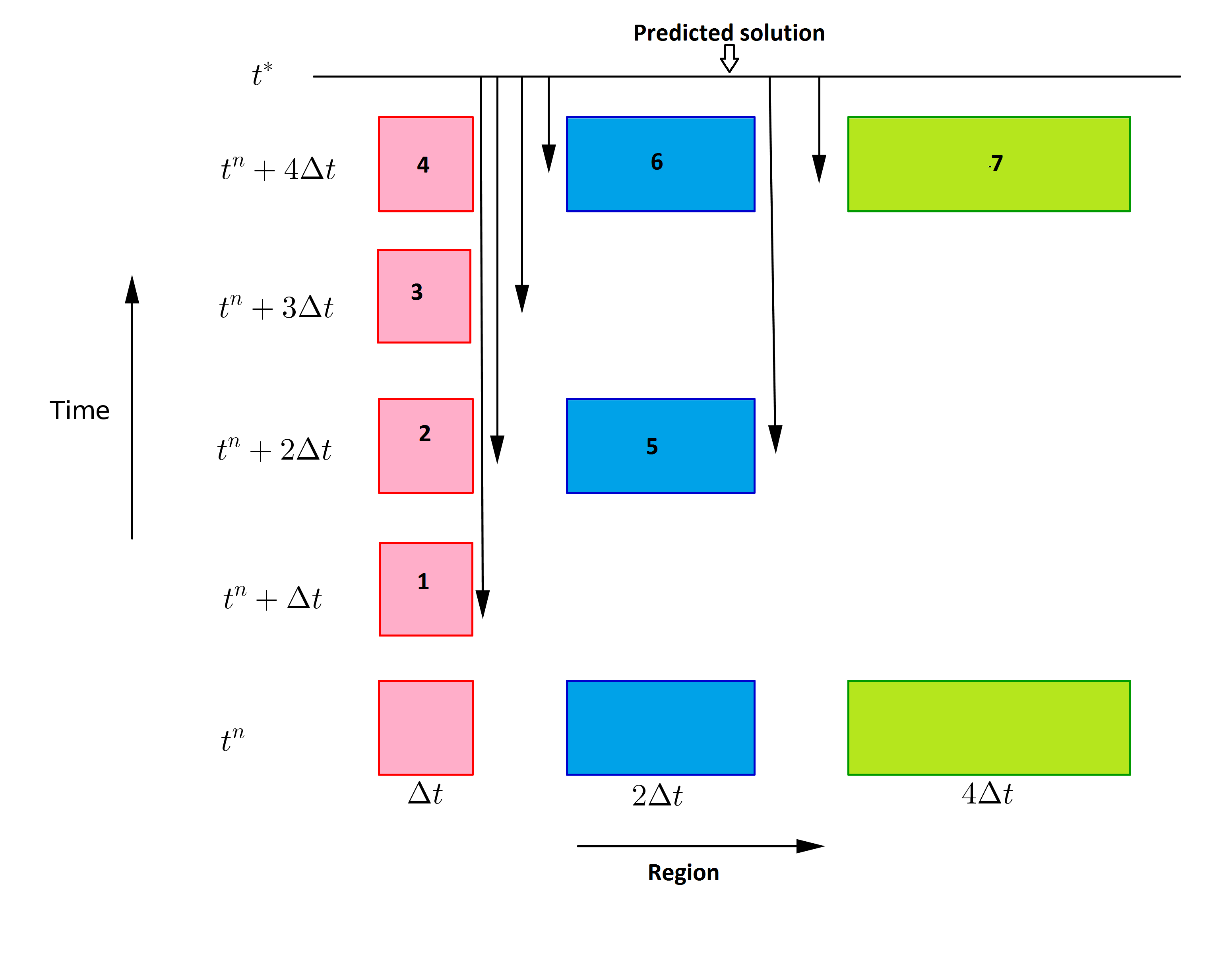}}
\caption{ (a) Prediction step and (b) correction step of the non overlap LTS method, when domain $\Omega$ is split into three regions $\Omega_0, \Omega_1$ and $\Omega_2$ respectively with time step $\Delta t, 2\Delta t$ and $4\Delta t$. }
\label{Non_over00}
\end{figure}

\begin{itemize}
\item \textbf{First step (Prediction):} advance the solution globally from  $t^n$ to $t^* = t^n + \Delta t^*$, by solving globally the DAREs using
the DG spatial discretization method and a time integrator with uniform time step $\Delta t^*$. This is schematically illustrated in \figref{\ref{Non_over00}} (a).

\item \textbf{Second step (Correction):} For all sub-domains $\Omega_i$, use the known component $X$ of the concentration at the time $t^n$ and $t^*$ to interpolate the value of $\Gamma_i$ at every time $t_i^j$, in order to advance the local component $X_i$ from time $t^n$ to $t^{n+1}$. This is schematically illustrated in \figref{\ref{Non_over00}} (b).
\end{itemize}
\newpage

This process can be repeated, in order to estimate the solution at any
final time $T^1$ from the known solution at  initial time $T^0$. To
that end, one can implement \Algref{\ref{NOver_pseudo_algo}} with $m_s =
\frac{T^1-T^0}{\Delta t_{\max}}$ and $r_i= \frac{\Delta t_{\max}}{\Delta t_i}$ for all $i\in\lbrace0,\cdots,m\rbrace$. Here, the function $\mathbf{I}$ is the interpolation function while the functions $\mathbf{A}_G$ and $\mathbf{A}_i$ are  the iterative function of the standard time integrators respectively used to solve the system of ODEs globally on $\Omega$ and locally on $\Omega_i$.  

\begin{algorithm}[H]
\caption{Pseudo algorithm of the non overlap LTS method.\label{NOver_pseudo_algo}}  
\begin{algorithmic}
\For{$i= 0 \text { to }  m-1$}
	\State $r_i=\frac{\Delta t_{\max}}{\Delta t_{i}}$\Comment{Number of steps needed to advance locally $X_i$ to a synchronized time using time step $\Delta t_i$}
\EndFor
\For{$n= 0 \text { to }  m^s-1$ }\Comment{Loop to advance the solution from $T^0$ to $T^{1}$}
    \State $t^n = T^0+n\times \Delta t_{\max}$\Comment{Compute the synchronized time $t^{n}$}
    \State $t^* = t^n+ \Delta t^*$\Comment{Compute the prediction time $t^{*}$}
    \State $X \Bigr|_{t^*} = \mathbf{A}_G\left(X \Bigr|_{t^{n}} \right) $\Comment{Computation of predicted global solution at time $t^*$}
    \For{$i= 0 \text { to }  m-1$}\Comment{Loop over all local regions - Correction step}
    	\State $t^{old}_i = t^n$\Comment{initialize the time of the known local solution $X_i$}
    	\For{$j= 1 \text { to }  r_i$ }\Comment{Loop to advance the solution  from $t^n$ to $t^{n+1}$ locally on $\Omega_i$}
    		\State $t^{j}_i = t^{old}_i + \Delta t_i$\Comment{Compute the local advancing time $t^{j}_i$}
    		\State $X\Bigr|_{\Gamma_i,t^{j}_i} = \mathbf{I}\left(  X \Bigr|_{ t^*}, X \Bigr|_{ t^n}\right) $\Comment{Interpolate the local solution $X_i$ on the boundary $\Gamma_i$ at time  $t_{new}^i$}
    		\State $X_i \Bigr|_{t^{j}_i} = \mathbf{A}_i\left(X_i \Bigr|_{t^{old}_i},  X\Bigr|_{\Gamma_i,t^{j}_i}  \right) $\Comment{Advance locally the solution on $\Omega_i$ from $t^{old}_i$ to $t^{j}_i$} 
    		\State $ t^{old}_i = t_{i}^j$\Comment{Update the time of the known local solution $X_i$}  
    	\EndFor
    \EndFor
\EndFor
\end{algorithmic}
\end{algorithm}

\section{Results and Discussion\label{ALL_num_res}}
The goal of this section is to investigate numerically the convergence
of LTS-DG schemes and compare their efficiency against a global time
stepping (GTS-DG) scheme. 
Firstly in \subsectref{\ref{OGAT_expe}}, by applying the OLTS-DG scheme to the two dimensional Ogata and Banks problem \cite{ogata}, we examine how the direction of the bulk velocity of the DAREs and the size of the overlap affect the accuracy of the  OLTS-DG schemes.
Secondly in \subsectref{\ref{OETO_expe}}, we compare the efficiency of the GTS-DG,
OLTS-DG and NOLTS-DG schemes, by applying them to the one dimensional
electron transfer only (ETO) model \cite{ET}. 
Finally, in \subsectref{\ref{TWO2}}, we examine the convergence and compare the efficiency of the GTS-DG and OLTS-DG when applied to the transport of solute through a 2D domain with fracture.

\subsection{Effect of the bulk velocity and the size of overlap on the OLTS-DG schemes\label{OGAT_expe}}

The purpose of this section is to investigate how the direction of the bulk velocity or the size of the overlap and the order in which the solution restraints to the  eligible sub-domains $S^r$ are consecutively solved, affect the accuracy of the numerical solution obtained with OLTS-DG schemes. To that end, the $DG_{\text{OLTSD-Impl}}$ and $DG_{\text{OLTSI-Impl}}$ schemes are used to solve the Ogata Banks equation with the bulk velocity $\beta = (1,0)$ and the diffusion coefficient $\epsilon = \frac{\parallel \beta\parallel}{Pe}$ where $Pe$ is the P\'eclet number.

\subsubsection{Effect of the bulk velocity on the OLTS-DG schemes}

In the Ogata and Banks problem \cite{ogata}, the fast change of the concentration of the solute takes place in a region close to the boundary at $x=0$. So the sub-domain that contains the boundary at $x=0$ should have the finest local time step, for a high accuracy of the OLTS-DG methods. This is illustrated by the better accuracy of both OLTS-DG methods obtained in the case where the sub-domain containing the boundary $x=0$ has the finest time step compared to case where it has the coarser time step (see \cite{AS17} for more details). Thus, to improve the efficiency of the OLTS-DG method, the choice of the fine, coarse discretized sub-domains  and the order of update of the local solution should respect the a priori physics.

\subsubsection{Effect of the size of overlap on the OLTS-DG schemes}
In this section, we investigate how the size of the overlap between
two sub-domains affect the accuracy of the global solution. To that
end, we consider the Ogata and Banks problem \cite{ogata} with the
initial sub-domains $\Omega_0 = [0,x_1]\times [0,1]$ and $\Omega_1 =
[x_1, x_2]\times [0,1]$.
For a given $n\in\{1,2,\cdots,11 \}$ and $h_x = 0.02$, we consider the overlapped sub-domains $\Omega_{0,n}$ and $\Omega_{1,n}$ given by
\begin{equation}
\Omega_{0,n} = [0,x_1+n\times h_x]\times [0,1], \qquad 
\Omega_{1,n} = [x_1-n\times h_x, x_2]\times [0,1].
\end{equation}
Note that the size of the overlap (i.e. $(\Omega_{0,n}\bigcap \Omega_{1,n}$) is equal to $2n\times h_x$  and increases with $n$. The local time steps are $\Delta t_0 = 2^{-11}$ and $\Delta t_1 = 2^{-10}$, thus we consider the $DG_{\text{OLTSI-Impl}}$ scheme for more accuracy.

\begin{figure}[H]
\centering
\includegraphics[height=5cm,width=10cm]{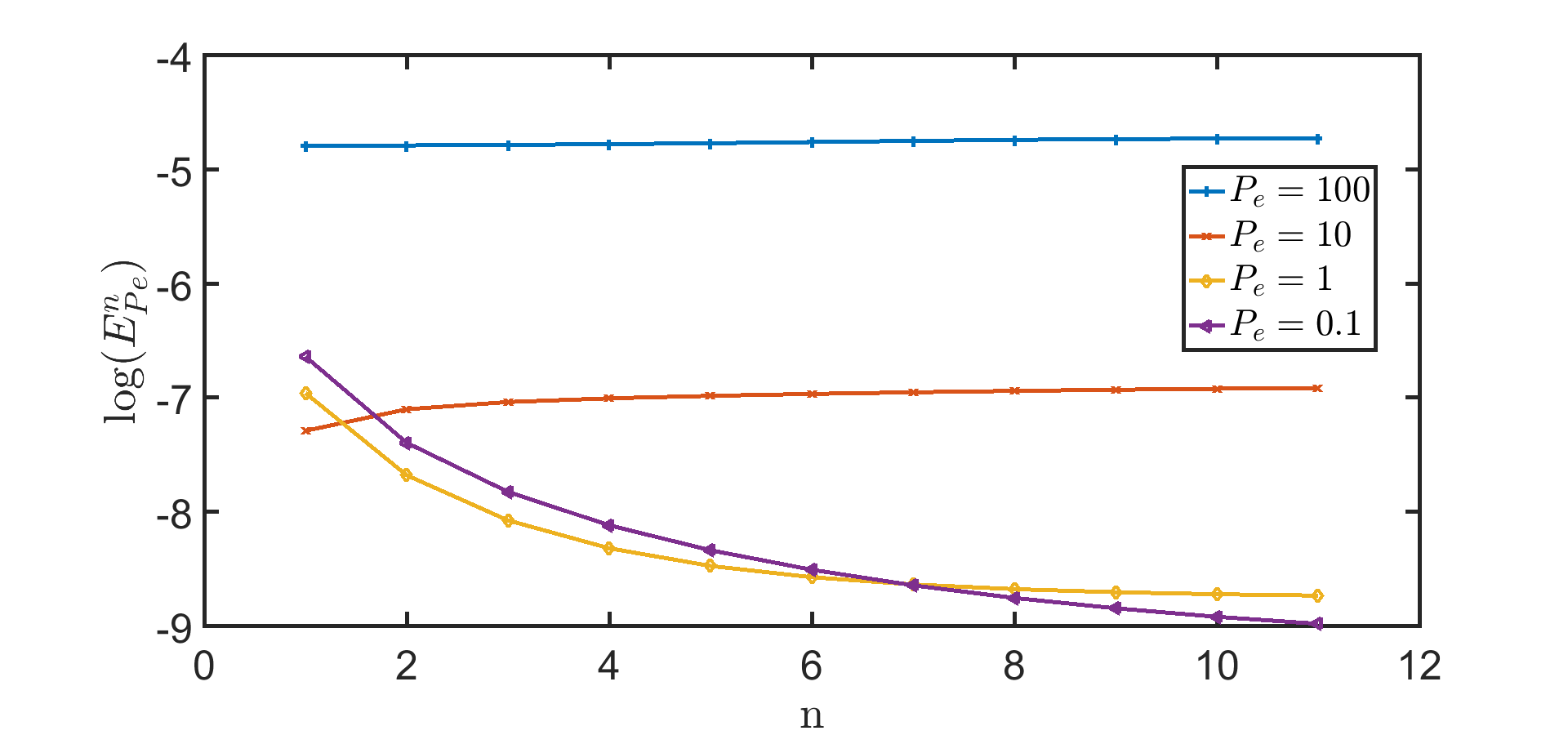}
\caption{We plot the logarithm relative error ($\log(E_{P_e}^n)$) against the integer $n$ for all P\'eclet number considered. The size of the overlap is equal to $2n\times h_x$ with $h_x = 0.02$}
\label{OvereffectFigRes}
\end{figure}

For all P\'eclet number $P_e \in \{0.1, 1, 10, 100 \}$ and all $n\in \{1,2,\cdots, 11\}$, we simulate the global solution,  $C_{P_e}^n$ at the time $t=0.5$, using the $DG_{\text{OLTSI-Impl}}$ scheme on the overlapped sub-domains $(\Omega_{0,n},\; \Delta t_0)$ and $(\Omega_{1,n},\; \Delta t_1)$. We then compute the relative error $E_{P_e}^n$ as follows
\begin{equation}
E_{P_e}^n = \frac{\Bigr| C_{P_e}^n - C \Bigr|_{L^2(\Omega)}}{\Bigr| C \Bigr|_{L^2(\Omega)}},
\end{equation}
where $C$ is the exact solution at the time $t=0.5$ and $\Omega$ the solution domain. The results are illustrated in \figref{\ref{OvereffectFigRes}}, where we plot the logarithm relative error ($\log(E_{P_e}^n)$) against the integer $n$ for all the P\'eclet numbers considered.
A few conclusions can be drawn from the results \figref{\ref{OvereffectFigRes}}: 
\begin{itemize}
\item For the P\'eclet number $P_e>1$ ($\epsilon<1$), the relative error $E_{P_e}^n$ increases slightly  as we increase the size of the overlap (i.e. as we increase $n$). Thus, in the case of high P\'eclet number, the overlapped sub-domains obtained by including only the direct neighbour into the initial sub-domains, is the best choice to simulate efficiently the global solution.
 
\item For the P\'eclet number $P_e\leq 1$ ($\epsilon\geq 1$), the relative error $E_{P_e}^n$ decreases as we increase the size of the overlap (i.e $n$).

\item For large size of overlap, i.e. $n>>1$, the relative error $E_{P_e}^n$ decreases as we decrease the P\'eclet number $P_e$ (or increase $\epsilon$). The same behaviour is observed in \cite{gander2005overlapping}, when the overlapping Schwarz waveform relaxation scheme was applied to the viscous Burger equation with various values of the viscosity parameter. For the overlap equal $0.2$  ($n=10$ in our case), the error decreases as the diffusion term increases.
\end{itemize}
\newpage

\subsection{Comparison of GTS-DG and LTS-DG schemes when applied to the 1D ETO model\label{OETO_expe}}
In this section, we compare the performance of the GTS-DG and LTS-DG schemes, when solving the one dimension ETO model. To that end,  we first describe how we split the global solution domain. Secondly, we investigate the derivation of the local system $\text{SO}_i$ associated to the sub-domain $\Omega_i$ for all $i=0,\cdots,m$. Finally, we present the numerical results. Here, we consider the LTS-DG schemes that use the Impl as time integrator. Thus we respectively denote $DG_{\text{OLTS-Impl}}$ and $DG_{\text{NOLTS-Impl}}$ the overlap and non overlap LTS-DG schemes.

Let us consider for example the ETO  process at the electrode represents by the chemical reaction defined as follow
\begin{equation}\label{ETO_reac}
\mathbf{Q} \xrightleftharpoons[k_b]{k_f} \mathbf{Q}^{+} +e^-,
\end{equation}
where the rate constants $k_b$ and $k_f$ are given  by the  Buttle-Volmer kinetics Equation \cite{ET}. The mathematical model to describe  ETO is derived from the Fick's Law for mass transfer \cite{Fick1, Fick2, Fick3}. In one dimension, the dimensionless governing equation  of ETO model is a coupled system of PDEs given by
\begin{subequations}
\begin{align}
\partial_{t}C_{\mathbf{Q}} =& \partial_z(\partial_zC_{\mathbf{Q}} ),\;\; &\forall(t,z)\in [0,2t_{\lambda}]\times[0,z_{\max}],\label{diffusion}\\
\partial_{t}C_{\mathbf{Q}^{+}} =& \partial_z(D^+\partial_zC_{\mathbf{Q}^{+}} ),\;\; &\forall(t,z)\in [0,2t_{\lambda}]\times[0,z_{\max}],\label{diffusion1}
\end{align}
\end{subequations}
subject to the boundary conditions
\begin{subequations}
\begin{align}
\partial_zC_{i}(t, z_{\max}) =& 0,\;\;\forall t\in [0,2t_{\lambda}], \;\;\forall i\in \{\mathbf{Q},\mathbf{Q}^{+}\},\label{diff1}\\
\partial_zC_{\mathbf{Q}}(t, 0)  =& K_f(t)C_{\mathbf{Q}}(t, 0) - K_b(t)C_{\mathbf{Q}^{+}}(t, 0),\;\;\forall t\in [0,2t_{\lambda}],\label{diff2}\\
D^+\partial_zC_{\mathbf{Q}^{+}}(t, 0)   =& -K_f(t)C_{\mathbf{Q}}(t, 0) + K_b(t)C_{\mathbf{Q}^{+}}(t, 0),\;\;\forall t\in [0,2t_{\lambda}],\label{diff3}
\end{align}
\end{subequations}
and the initial condition
\begin{equation}\label{diff4}
C_{\mathbf{Q}}(0, z) =1,\;\;C_{\mathbf{Q}^{+}}(0, z) =0    ,\;\; \forall z\in [0,z_{\max}].
\end{equation}
Here $C_{i}(t,z)\in\mathbb{R}$ is the dimensionless concentration of the species $i\in \{\mathbf{Q},\mathbf{Q}^{+}\}$ and $D^+\in\mathbb{R}$ is the dimensionless diffusion of the specie $\mathbf{Q}^{+}$. The heterogeneous electron transfer rate constants, $k_b, k_f$, can  be written in its dimensionless form, $K_b, K_f$,  as follows
\begin{equation}\label{reduc}
K_f = K_0\exp[(1-\delta)P],\qquad
K_b = K_0\exp[(-\delta) P],
\end{equation}
where the dimensionless potential, $P$, in terms of the dimensionless time, is given by
\begin{align}
P = \left\lbrace
\begin{array}{ll}
P_1 + t , \; 0\leq t\leq t_\lambda\\
P_2 -  (t - t_\lambda) , \; t_\lambda\leq t\leq 2t_\lambda
\end{array} 
\right.,\qquad
t_\lambda = P_2-P_1,
\end{align}
with $P_1$ and $P_2$ respectively the dimensionless initial and reverse potential. The dimensionless current, $G$, is given by 
\begin{equation}\label{cur_dim}
G(t) =   \partial_zC_{\mathbf{Q}}(t, 0)  = K_f(t)C_{\mathbf{Q}}(t, 0) - K_b(t)C_{\mathbf{Q}^{+}}(t, 0),\;\;\forall t\in [0,2t_{\lambda}].
\end{equation}

The solution domain is given by $\Omega = [0,z_{\max}]$, with $z_{\max}$ proportional to the diffusion length $\delta$ (i.e. $z_{\max}= k\delta,\; k\in \mathbb{N})$. Since the dimensionless current depends only on the concentration of the species at the boundary $z=0$, we consider the partition $ \Omega = \cup_{i = 1}^{n}I_i$ where the interval $I_i=[z_{i-1},z_i]$ is such that the step size $h_i=z_i-z_{i-1}$, respectively, follows the geometric and the uniform progression on $[0,\delta]$ and $[\delta, z_{max}]$. Specifically for a given number $r\in \mathbb{N}$ and the increasing factor $q$, we have 
\begin{equation*}
h_i= 
\left\lbrace 
\begin{array}{ll}
h\times q^{i-1},\;& i= 1, \cdots, r\\
h\times q^{r},\;&i= r+1, \cdots, n
\end{array} 
\right.,
h = \delta\left( \frac{q^{r} - 1}{q-1}\right)^{-1}, n=r+\lceil\frac{z_{\max}-\delta}{h_{r}}\rceil. 
\end{equation*}
Unless stated, for the simulation we use $q=1.05$ and $r=100$. We associate to each element $I_i$ the time step $\Delta T_{I_i}\leq 2^{n_i}, \; n_i=\lceil\log_2\left(\mathbf{C}_{max}h_i \right)\rceil$ for a given Courant number $\mathbf{C}_{max}$. We then update the time step on each element $I_i$ by setting $\Delta T_{I_i} = \Delta T_{I_1}  $ for all $i=1,\cdots,r$ and $\Delta T_{I_i} = \Delta T_{I_{r+1}}$ for all $i=r+1,\cdots,n$. 

For the simulation, we consider the geometry settings and the ETO model parameters given by 
\begin{equation}\label{ETO_LTS_settings}
\delta = 20, z_{\max} = 5\delta, \mathbf{C}_{max} = 0.3, K_0=20, \tilde{D}^+=1.
\end{equation}
This leads to the local time steps $\Delta t_{I_1} \leq 2^{-9}$ and $\Delta t_{I_{r+1}} \leq 2^{-2}$. We finally consider the sub-domains $\Omega_0 = [0,z_r]$ and $\Omega_1=[z_r,z_{n}]$ with the local time step $\Delta t_0 = 2^{-9}$ and $\Delta t_1 = 2^{-6}$, respectively.

By considering an orthonormal DG finite space, the DG space discretization of the governing equation of the ETO model leads to the system  of ODEs defined as follows

\begin{equation}\label{syst_L}
\underbrace{\left( 
\begin{array}{c}
d_t\alpha\\
d_t\alpha^+
\end{array}
\right)}_{d_t\chi}
 + 
\underbrace{\left( 
\begin{array}{cc}
L^{s,\sigma_0} + K_f(\tilde{t})F^1 &-K_b(\tilde{t})F^1\\
-K_f(\tilde{t})F^1&D^+L^{s,\sigma_0}+K_b(\tilde{t})F^1
\end{array}
\right)}_{\mathcal{L}^{s,\sigma_0}_{D^+,K_f,K_b}}
\underbrace{\left( 
\begin{array}{c}
\alpha\\
\alpha^+
\end{array}
\right)}_{\chi} =
\left( 
\begin{array}{c}
0\\
0
\end{array}
\right).
\end{equation}
Here $\alpha$ and $\alpha^+$ are respectively the components of the concentration of the species $Q$ and $Q^+$ in the DG finite space. The matrix $L^{s,\sigma_0}$ represents the DG driscretization of the diffusion operator. More details on this derivation can be found in \cite{AS17}.

\subsubsection{Local ODE system for the overlap LTS-DG scheme}
\begin{figure}[H]
\centering
\includegraphics[height=6cm,width=10cm]{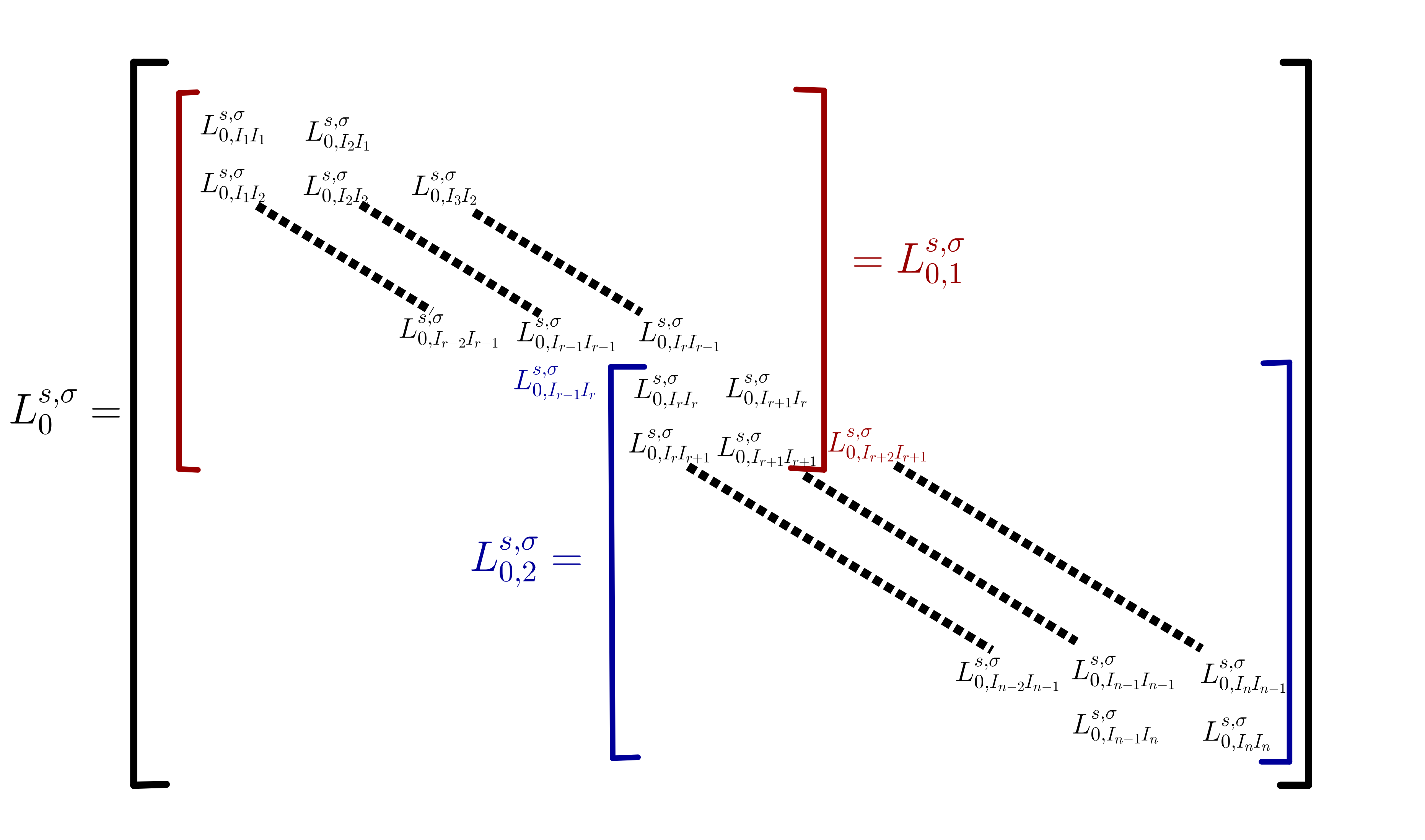}
\caption{Extraction of the matrices $L_{0,1}^{s,\sigma_0}$ and $L_{0,1}^{s,\sigma_0}$ from the stiffness matrix $L_{0}^{s,\sigma_0}$ when applying the overlap LTS method to the one dimension ETO model. It shows that only the block matrices $  L_{0,I_{r+2}I_{r+1}}^{s,\sigma} $ and $  L_{0,I_{r-1}I_{r}}^{s,\sigma} $, from the matrix $L_{0}^{s,\sigma_0}$ respectively contribute to the computation of the vectors $\mathbf{B}_1^T$ and $\mathbf{B}_2^T$.}
\label{over_mat_extrat}
\end{figure}

In this section, we show how to extract the local ODE system for the
OLTS-DG scheme from the global ODE system of the 1D ETO model. To
overlap the sub-domains here, we include the direct neighbour into the
initial sub-domains. Thus, the overlapped sub-domains are $\Omega_0 =
[0, z_{r+1}]$ and $\Omega_1 = [ z_{r-1}, z_n]$.
In this case, we consider as internal boundary $\Gamma_0$ and $\Gamma_1$ as the node $z_{r+1}$ of $[z_{r+1}, z_{r+2}]$ and $z_{r-1}$ of $[z_{r-2}, z_{r-1}]$, respectively.
The system of ODEs given by \Equaref{\ref{syst_L}}, obtained from the DG spatial discretization of the dimensionless governing equation of the ETO model, can be split into two systems of ODEs 
\begin{align}
\dfrac{d \chi_0}{d\tilde{t}} = L_0(\tilde{t}) \chi_0 + \mathbb{S}_0^e\chi_1\Bigr|_{\Gamma_0},\label{ETO_ODES_sp1} \\
\dfrac{d \chi_1}{d\tilde{t}} = L_1(\tilde{t}) \chi_1 + \mathbb{S}_1^e\chi_0\Bigr|_{\Gamma_1},\label{ETO_ODES_sp2}
\end{align}
where $\chi_j \in \mathbb{R}^{2n_{\Omega_j}}$ is the coupled component of the concentration of the species $\mathbf{Q}, \mathbf{Q}^+$ on the region $\Omega_j$ for all $j=0,1$. The dimensions $n_{\Omega_0}$ and $n_{\Omega_1}$ are given by  
\begin{equation}
n_{\Omega_0} = \sum_{i=1}^{r+1} (k_i+1), \qquad n_{\Omega_1} = \sum_{i=r}^n (k_i+1),
\end{equation}
where $k_i$ is the highest degree of the Legendre polynomials considered on $I_i$ (i.e. $k_i+1$ is the dimension of the DG finite space on $I_i$). The matrices $L_j, \mathbb{S}_j^e \in \mathbb{R}^{2n_{\Omega_j}\times 2n_{\Omega_j}}$ for all $j=0,1$ can be obtained from the matrix $\mathcal{L}^{s,\sigma_0}_{D^+,K_f,K_b}$,  given by \Equaref{\ref{syst_L}}, as follows
\begin{equation}\label{split_over_ETO_matr}
L_1=\left( 
\begin{array}{c|c}
DL_{0,0}^{s,\sigma_0} + K_f(\tilde{t})L^1 &-K_b(\tilde{t})L^1\\\hline
-K_f(\tilde{t})L^1&D^+L_{0,0}^{s,\sigma_0}+K_b(\tilde{t})L^1
\end{array}
\right), \;
L_2=\left( 
\begin{array}{c|c}
DL_{0,1}^{s,\sigma_0}&0\\\hline
0&D^+L_{0,1}^{s,\sigma_0}
\end{array}
\right).
\end{equation}
Here, the matrix $L^1\in \mathbb{R}^{n_{\Omega_0}\times n_{\Omega_0}}
$ takes the same form as the matrix $F^1$ defined in \Equaref{\ref{syst_L}}; the
matrices $L_{0,0}^{s,\sigma_0}$ and $L_{0,1}^{s,\sigma_0}$ are
efficiently extracted from the matrix $L_{0}^{s,\sigma_0}$ as
illustrated by \figref{\ref{over_mat_extrat}}, due to its
tridiagonalisation. Note from \figref{\ref{over_mat_extrat}} that only the
block matrices $  L_{0,I_{r+2}I_{r+1}}^{s,\sigma} $ and $
L_{0,I_{r-1}I_{r}}^{s,\sigma} $, from the matrix $L_{0}^{s,\sigma_0}$
respectively contribute to the computation of the vector
$\mathbf{B}_0^T =\mathbb{S}_0^e\chi_1\Bigr|_{\Gamma_0}$ and
$\mathbf{B}_1^T =\mathbb{S}_1^e\chi_0\Bigr|_{\Gamma_1}$.
\newpage
For all $j\in \lbrace 0,1\rbrace$, we have $\mathbf{B}_j = [\mathbf{B}_j^{\alpha}, \mathbf{B}_j^{\alpha^+}]$ where the transpose of the block vectors $\mathbf{B}_j^{\alpha}, \mathbf{B}_j^{\alpha^+}\in \mathbb{R}^{1\times n_{\Omega_j}}$ are given by 
\begin{equation}\label{inter_boun_cont}
\mathbf{B}_0^{\zeta} = \left(  \mathbf{B}_{0,I_i}^{\zeta} \right)_{i=1,\cdots,r+1}, \quad 
\mathbf{B}_1^{\zeta} = \left(  \mathbf{B}_{1,I_{p_i}}^{\zeta} \right)_{i=1,\cdots,n-r+1} ,\; p_i = i+r-1,
\end{equation}
for all $\zeta=\alpha, \alpha^+$. According to the extraction of matrices illustrated in \figref{\ref{over_mat_extrat}}, we have for all $\zeta=\alpha, \alpha^+$
\begin{align*}
\mathbf{B}_{0,I_i}^{\zeta} =& 
\left\lbrace 
\begin{array}{ll}
D^{\zeta}L_{0,I_{r+2}I_{r+1}}^{s,\sigma} \zeta_{I_{r+2}}   \;\;&\text{if}\;\;i=r+1\\
0\;\;&\text{if}\;\;i=1,\cdots, r
\end{array}
\right.,\\
\mathbf{B}_{1,I_{p_i}}^{\zeta} = &
\left\lbrace 
\begin{array}{ll}
D^{\zeta}L_{0,I_{r-1}I_{r}}^{s,\sigma}\zeta_{I_{r-1}}  \;\;&\text{if}\;\;i=1\\
0\;\;&\text{if}\;\;i=2,\cdots, n-r+2
\end{array}
\right.,
\end{align*}
where the coefficient $D^{\zeta}$ is the dimensionless diffusion coefficient of the species with the component in the DG finite space $\zeta=\alpha, \alpha^+$.

\subsubsection{Local ODE system for the non overlap LTS-DG scheme}
\begin{figure}[H]
\centering
\includegraphics[height=6cm,width=10cm]{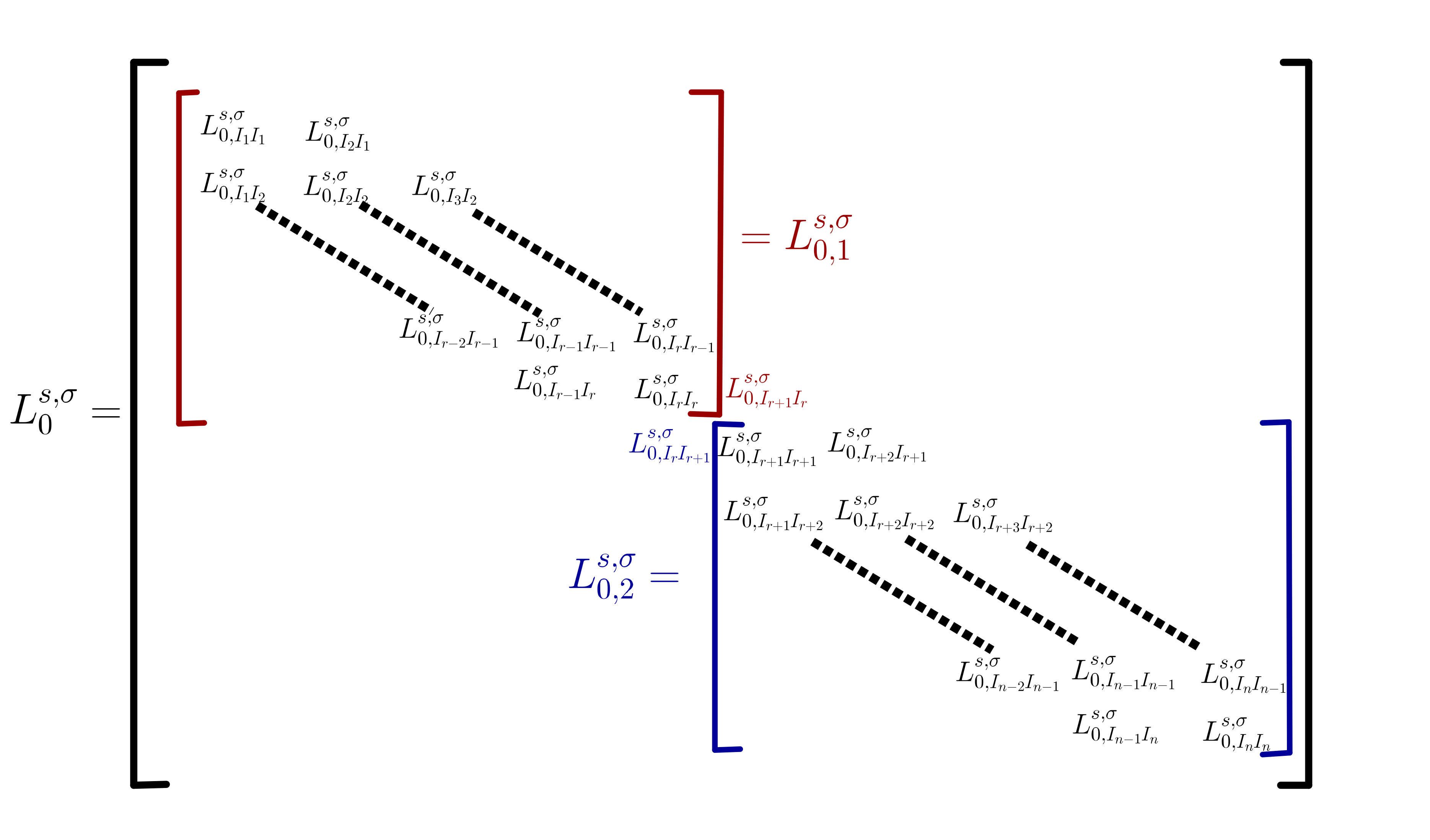}
\caption{Extraction of the matrices $L_{0,1}^{s,\sigma_0}$ and $L_{0,1}^{s,\sigma_0}$ from the matrix $L_{0}^{s,\sigma_0}$ when applying the non overlap LTS method to the one dimension ETO model. It shows that only the block matrices $  L_{0,I_{r+1}I_{r}}^{s,\sigma} $ and $  L_{0,I_{r}I_{r+1}}^{s,\sigma} $, from the matrix $L_{0}^{s,\sigma_0}$ respectively contribute to the computation of the vectors $\mathbf{B}_1^T$ and $\mathbf{B}_2^T$.}
\label{nover_mat_extrat}
\end{figure}

In this section, we show how to extract the local ODE system for the
NOLTS-DG scheme from the global ODE system of the 1D ETO model. The
non overlapped sub-domains are $\Omega_0 = [0, z_{r}]$ and $\Omega_1 =
[ z_{r}, z_n]$.
In this case, we consider the node $z_r$ of the interval $[z_{r}, z_{r+1}]$ and $[z_{r-1}, z_{r}]$ as the internal boundary of $\Gamma_0$ and $\Gamma_1$ , respectively.
The system of ODEs obtained from the DG spatial discretization of the dimensionless governing equation of the ETO model, can be split into two ODEs system of \Equaref{\ref{ETO_ODES_sp1}} and \Equaref{\ref{ETO_ODES_sp2}}. In this case, the dimensions $n_{\Omega_0}$, $n_{\Omega_1}$ are given by  
\begin{equation}
n_{\Omega_0} = \sum_{i=1}^{r} (k_i+1), \qquad n_{\Omega_1} = \sum_{i=r+1}^n (k_i+1),
\end{equation}
where $k_i$ is the highest degree of the Legendre polynomials considered on $I_i$. Also, the matrices $L_j, \mathbb{S}_j^e \in \mathbb{R}^{2n_{\Omega_j}\times 2n_{\Omega_j}}$ for all $j=0,1$ are given by \Equaref{\ref{split_over_ETO_matr}}
where the matrix $L^1\in \mathbb{R}^{n_{\Omega_0}\times n_{\Omega_0}} $ takes the same form as the matrix $F^1$ defined in \Equaref{\ref{syst_L}}; the matrices $L_{0,0}^{s,\sigma_0}$ and $L_{0,1}^{s,\sigma_0}$ are efficiently extracted from the matrix $L_{0}^{s,\sigma_0}$ as illustrated by \figref{\ref{nover_mat_extrat}}, due to its tridiagonalisation. 
\newpage
Note from \figref{\ref{nover_mat_extrat}} that only the block matrices $  L_{0,I_{r+1}I_{r}}^{s,\sigma} $ and $  L_{0,I_{r}I_{r+1}}^{s,\sigma} $, from the matrix $L_{0}^{s,\sigma_0}$ respectively contribute to the computation of the vector $\mathbf{B}_0^T =\mathbb{S}_0^e\chi_1\Bigr|_{\Gamma_0}$ and $\mathbf{B}_1^T =\mathbb{S}_1^e\chi_0\Bigr|_{\Gamma_1}$. For all $j\in \lbrace 0,1\rbrace$, we have $\mathbf{B}_j = [\mathbf{B}_j^{\alpha}, \mathbf{B}_j^{\alpha^+}]$ where the transposes of the block vectors $\mathbf{B}_j^{\alpha}, \mathbf{B}_j^{\alpha^+}\in \mathbb{R}^{1\times n_{\Omega_j}}$ are given by 
\begin{equation}
\mathbf{B}_0^{\zeta} = \left(  \mathbf{B}_{0,I_i}^{\zeta} \right)_{i=1,\cdots,r}, \quad 
\mathbf{B}_1^{\zeta} = \left(  \mathbf{B}_{1,I_{p_i}}^{\zeta} \right)_{i=1,\cdots,n-r} ,\; p_i = i+r.
\end{equation}
According to the extraction of matrices illustrated in \figref{\ref{nover_mat_extrat}}, we have for all $\zeta=\alpha, \alpha^+$,
\begin{align*}
\mathbf{B}_{0,I_i}^{\zeta} =& 
\left\lbrace 
\begin{array}{ll}
D^{\zeta}L_{0,I_{r+1}I_{r}}^{s,\sigma}\zeta_{I_{r+1}} \;\;&\text{if}\;\;i=r\\
0\;\;&\text{if}\;\;i=1,\cdots, r-1
\end{array}
\right.,\\
\mathbf{B}_{1,I_{p_i}}^{\zeta} = &
\left\lbrace 
\begin{array}{ll}
D^{\zeta}L_{0,I_{r}I_{r+1}}^{s,\sigma}\zeta_{I_{r}}  \;\;&\text{if}\;\;i=1\\
0\;\;&\text{if}\;\;i=2,\cdots, n-r+1
\end{array}
\right..
\end{align*}

\subsubsection{Numerical results of GTS-DG and LTS-DG schemes applied to the ETO model}
Let us now focus on the numerical comparison of the accuracy and the efficiency of the GTS-DG and LTS-DG schemes, while simulating the dimensional current of the ETO model. To that end, for a given $i\in \{0,\cdots,4\}$, we simulate the dimensionless current $G^{i,q},$ $q=DG_{\text{Impl}}$, $DG_{\text{OLTS-Impl}}$, $DG_{\text{NOLTS-Impl}}$ for the local time step $\Delta t_j^{i} = 2^{-i}\times \Delta t_j $ on the local solution domain $\Omega_j$ for all $j=0,1$. Note that for a given $i\in \{0,\cdots,4\}$, the universal time step, $ht^i=\Delta t_0^{i}$ , is considered for the GTS-DG schemes (i.e. the finest local time step of LTS-DG schemes). During the simulation of $G^{i,q}$, we also record the computation time $\text{CPU}^{i,q}$.

\begin{figure}[H]
\centering
\subfloat[]{\includegraphics[height=0.2\textheight,width=0.48\textwidth]{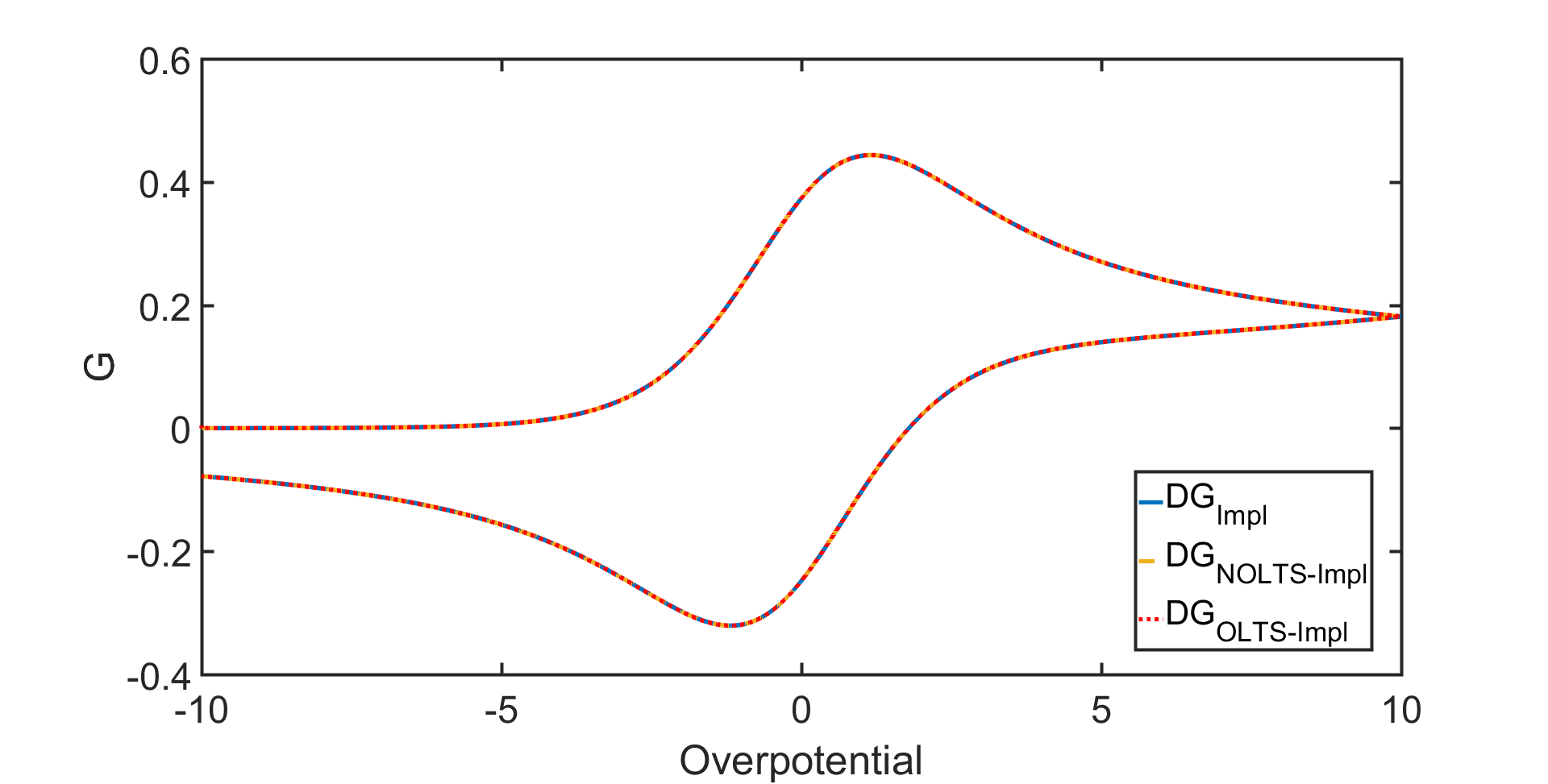}}
\subfloat[]{\includegraphics[height=0.2\textheight,width=0.48\textwidth]{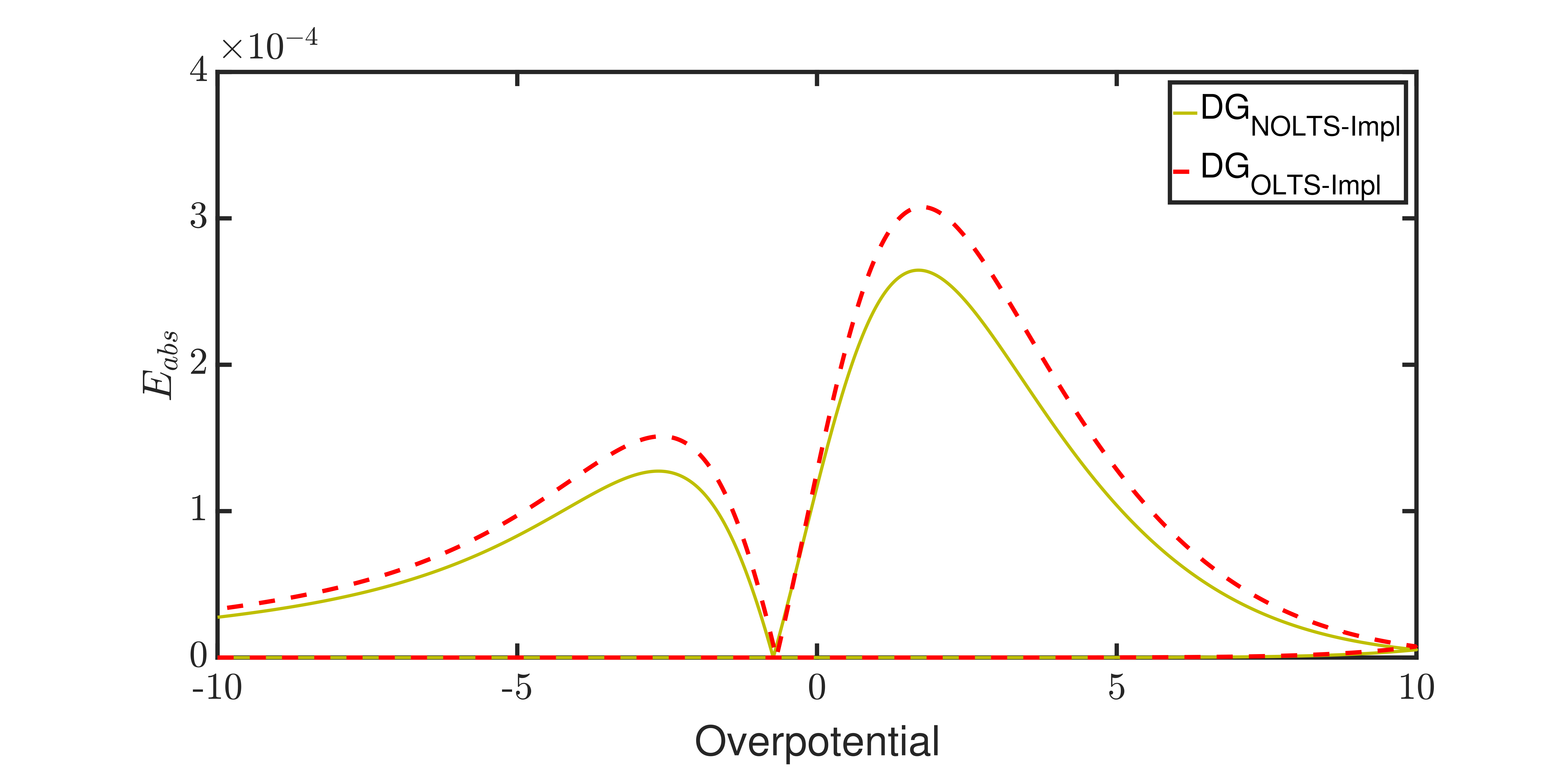}}   
\caption{\textbf{Voltammogram of the ETO model simulate with GTS-DG and LTS-DG schemes.} In (a), we plot $G^{0,q}$ against the overpotential for all $q=DG_{\text{Impl}}, DG_{\text{OLTS-Impl}}, , DG_{\text{NOLTS-Impl}}$. In (b), we plot the absolute  difference $E_{\text{abs}}^q$ against the overpotential for all $q=DG_{\text{OLTS-Impl}}, , DG_{\text{NOLTS-Impl}}$.
\label{overlap_LTS_solution_ETO}}
\end{figure}

By assuming that the exact dimensionless current is given by $G^{4,DG_{\text{Impl}}}$, we then compare it against $G^{4,q}$ $q= DG_{\text{OLTS-Impl}}$, $DG_{\text{NOLTS-Impl}}$ to see which of the LTS-DG schemes is more accurate. This is illustrated in \figref{\ref{overlap_LTS_solution_ETO}}, by plotting the dimensionless currents $G^{4,q}$ against the overpotential for all $q=DG_{\text{Impl}}$, $DG_{\text{OLTS-Impl}}$, $DG_{\text{NOLTS-Impl}}$ in \figref{\ref{overlap_LTS_solution_ETO}}(a); and the absolute  difference $E_{\text{abs}}^q = \Bigr| G^{4,DG_{\text{Impl}}} - G^{4,q} \Bigr|$ against the overpotential for all $q= DG_{\text{OLTS-Impl}}, DG_{\text{NOLTS-Impl}}$ in \figref{\ref{overlap_LTS_solution_ETO}}(b).

\newpage

\begin{figure}[H]
\centering
 \subfloat[]{\includegraphics[height=0.2\textheight,width=0.48\textwidth]{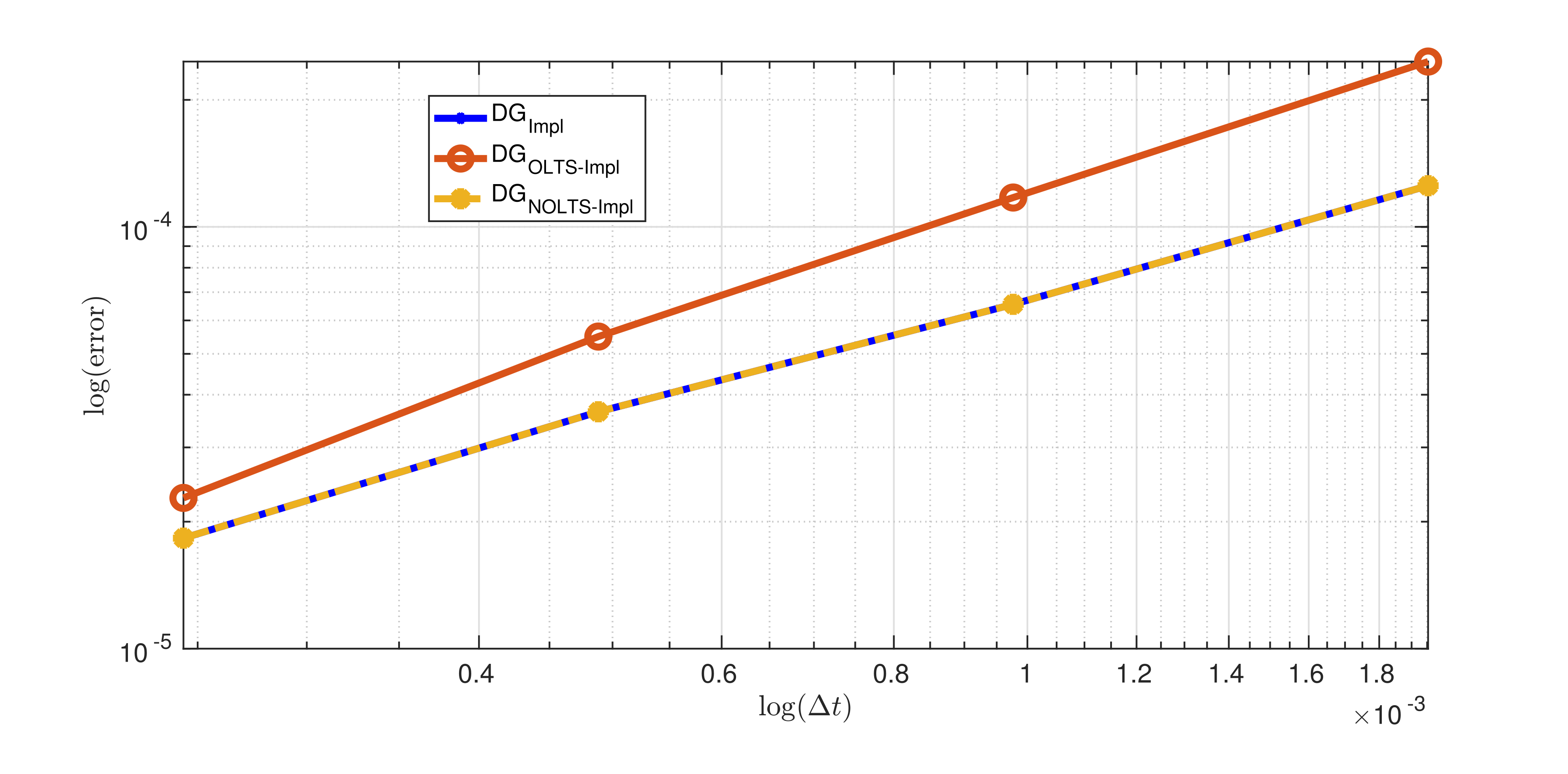}}
\subfloat[]{\includegraphics[height=0.2\textheight,width=0.48\textwidth]{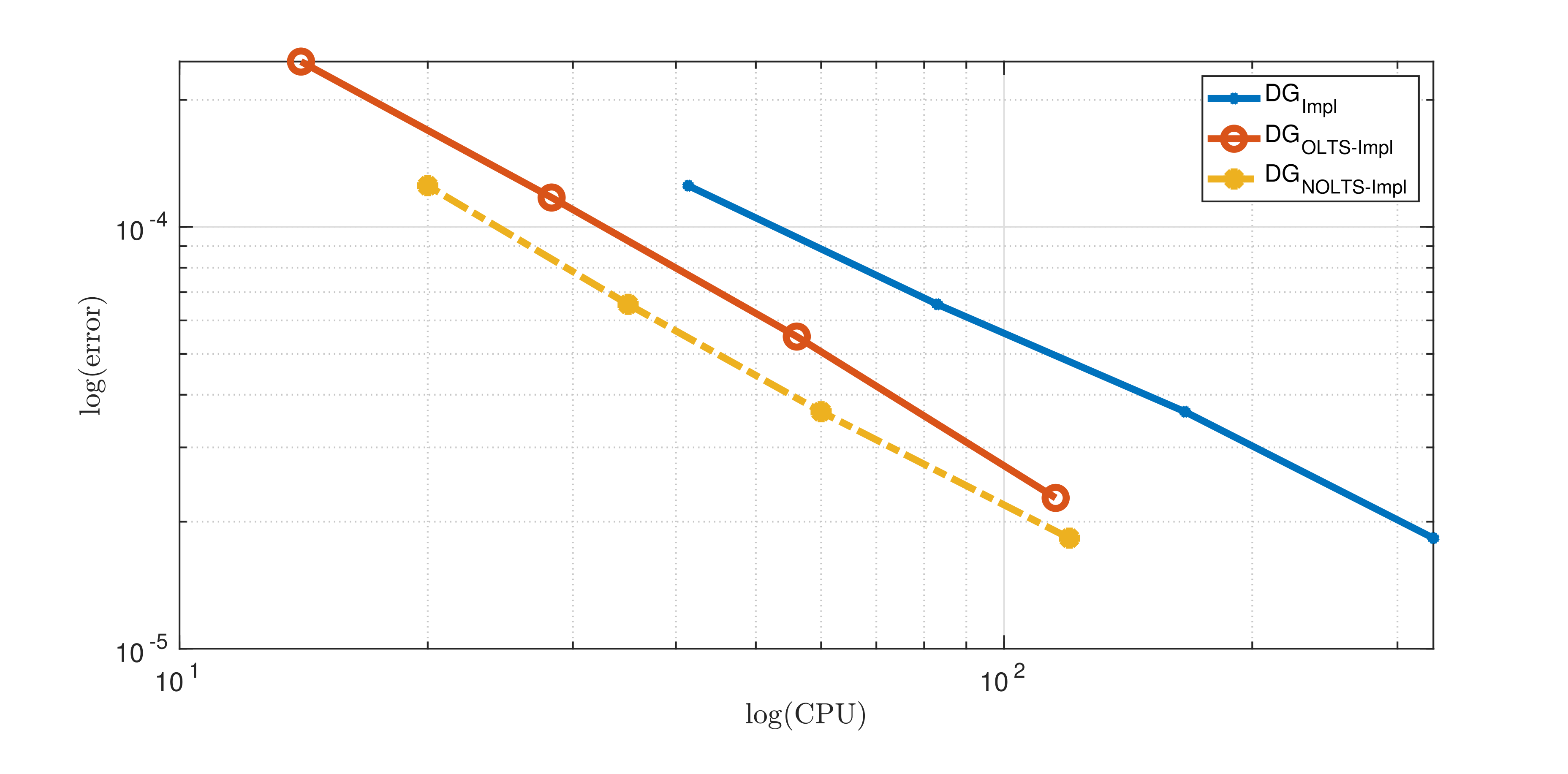}}  
\caption{Convergence and efficiency of the GTS-DG and LTS-DG schemes
  for the ETO model. We plot in (a) the error
  $\log(\text{error}^{i}_{q})$ against the minimum time step
  $\log(ht^{i})$ and in (b) the computation time $\log(\text{CPU}^{i,q})$ for all $i=0,\cdots,3$ and all $q=DG_{\text{Impl}}, DG_{\text{OLTS-Impl}}, DG_{\text{NOLTS-Impl}}$. 
\label{overlap_LTS_perfor_ETO}}
\end{figure}

To investigate the convergence and the efficiency of the GTS-DG and LTS-DG schemes, we compute the relative errors, $\text{error}^{i}_{q}$ , given by
\begin{equation}
\text{error}^{i}_{q} = 100\dfrac{ \Bigr| G^{4,DG_{\text{Impl}}} -  G^{i,q}\Bigr|_{L^2(P)}^2  }{ \Bigr| G^{4,DG_{\text{Impl}}} \Bigr|_{L^2(P)}^2},
\end{equation}
for all $i=0,\cdots,3$ and all solver $q=DG_{\text{Impl}}, DG_{\text{OLTS-Impl}}, DG_{\text{NOLTS-Impl}}$.    
We then plot in \figref{\ref{overlap_LTS_perfor_ETO}}(a) the error , $\log(\text{error}^{i}_{q})$, against the minimum of the local time step, $\log(ht^{i})$. This shows the decay of the error with respect to the minimum local time step, meaning the  $DG_{\text{Impl}}$,  $DG_{\text{OLTS-Impl}}$, $DG_{\text{NOLTS-Impl}}$ converge in time. In \figref{\ref{overlap_LTS_perfor_ETO}}(b), we plot
the error, $\log(\text{error}^{i}_{q})$, against the computation time, $\log(\text{CPU}^{i,q})$.  Note from \figref{\ref{overlap_LTS_perfor_ETO}}(b) that 
for a given error $E\in \mathbb{R}$ such that $E = \text{error}^i_q$, $q=DG_{\text{Impl}}, DG_{\text{OLTS-Impl}}, DG_{\text{NOLTS-Impl}}$, we have
$$\text{CPU}^{i, DG_{\text{NOLTS-Impl}}} <\text{CPU}^{i, DG_{\text{OLTS-Impl}}} < \text{CPU}^{i, DG_{\text{Impl}}}.$$ 
\figref{\ref{overlap_LTS_perfor_ETO}}(b) shows that the LTS-DG schemes, compared to GTS-DG schemes, are more efficient to simulate the dimensionless current of the  ETO model.

The 1D numerical experiment realized in this section has shown that 
\begin{itemize}
\item even if the OLTS-DG schemes is not as accurate as GTS-DG scheme, the computation time of the OLTS-DG schemes is small enough to make them more efficient compare to the GTS-DG schemes,
\item the non overlap LTS-DG scheme is more accurate than the overlap LTS-DG scheme. This is because the estimation of the needed internal value $ \Gamma_i $, for the resolution of the system $\text{SO}_i$, is more accurate with the non overlap LTS-DG scheme.
\end{itemize}

\subsection{Comparison of GTS-DG and OLTS-DG schemes when applied to the 2D transport of solute through a domain with fracture\label{TWO2}}

Consider the transport of an inert solute within an incompressible
fluid with an absence of volumetric source and sinks, through a 2D
domain, $\Omega$ with a fracture.
The concentration $C(x,y,t)$ of the solute follows 

\begin{equation}\label{RDA51}
\dfrac{\partial C}{\partial t}  - \nabla \cdot (D_m\nabla C)+ \nabla\cdot\textbf{v}C = R(C) \;\;\text{ with }\;\; R(C) = 0,
\end{equation}
where $R(C)$ is the reaction term, $D_m$ is a molecular diffusivity and the velocity $\textbf{v}$ in each pore, computed from the solution of Darcy's equation with the equation pressure given by
\begin{equation}\label{pressure_eq}
\textbf{v} =  -\frac{\textbf{k}}{\mu\phi}\nabla p
\;\;\text{ with }\;\;\left\lbrace \begin{array}{lll}
\nabla \cdot \left(-\frac{k}{\mu}\nabla p \right) &=0, &\\
p &= p_0 \quad &\text{in} \quad \partial \Omega_D^1,\\
\vec{n}\cdot\left(-\frac{k}{\mu}\nabla p \right)  &= p_1 \quad &\text{in} \quad \partial \Omega_N^1. 
\end{array}
\right.
\end{equation}
Here $\textbf{v}$ is the fluid velocity, $p$ is the pressure, $\phi$ is the porosity, $k$ is the permeability, $\mu$ is the viscosity.
As a boundary and initial condition, we keep the concentration $C$ and
the pressure $p$ respectively at a constant value $C_0$ and $p_0$ at
the inflow boundary $\partial\Omega_2$ and allow it to undergo pure
advection at the outflow boundary $\partial\Omega_4$. The boundary
conditions also include the no flux at the rigid boundaries,
The fracture is represented by the domain $\Omega_r = [x_r, x_r + l]\times[y_r, y_r-h]$. We assumed that within the fracture, the permeability is 1000 times greater than the permeability of the remaining domain. 

For the simulation of the concentration profile we  design the
unstructured mesh of the domain for $L=1$ with $distmesh2d$, such that
for a given element $T$ within the fracture we have the radius on the
incircle of $T$ is less than $h/3$. The finest elements are located in
the fracture characterized by the coordinates $x_r = 0.2, y_r = 0.51$,
the height $h=0.02$ and the length $l=0.6$. This is illustrated in
\figref{\ref{frac}} (a). 

\begin{figure}[H]
  \centering         
  \subfloat[]{\includegraphics[height=0.15\textheight,width=0.48\textwidth]{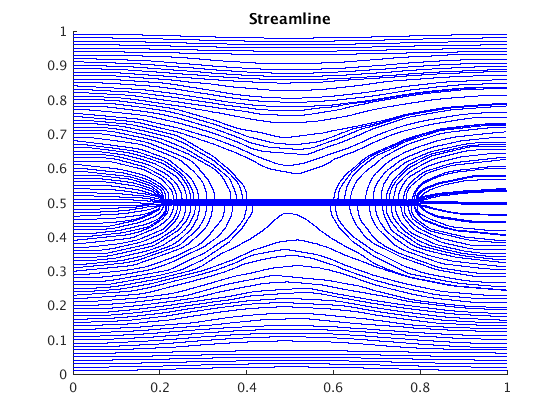}}
  \subfloat[]{\includegraphics[height=0.15\textheight,width=0.48\textwidth]{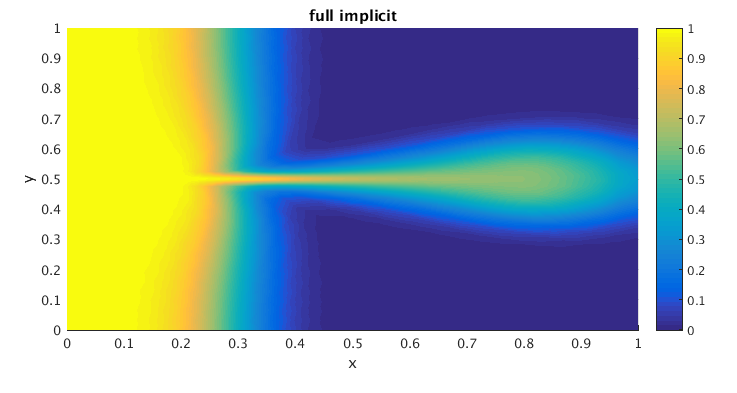}}
  \caption{(a) Streamline of the simulated fluid velocity of solute
    through domain with fracture. Note from this figure that the
    velocity of the fluid is higher within the fracture.  
    (b) Solution without presence of source or reaction ($R=0$) found
  using full implicit.}
\label{velo_concen3}
\end{figure}

In this case, the equation of the pressure is a steady heterogeneous
diffusion equation given by \Equaref{\ref{pressure_eq}}  with 
\begin{equation}
p_0 = 1,\; p_1 = 0,\; k\mu^{-1}=
\left\lbrace \begin{array}{ll}
1000 &\quad\text{on}\quad \Omega_r\\
1&\quad\text{on}\quad \Omega\setminus\Omega_r
\end{array}
\right. .
\end{equation}
We then simulate the fluid velocity on each element $T$ after using
the SIPG method to solve the equation of the pressure. 
For the sake of clarity, we plot in \figref{\ref{velo_concen3}}(a) the
streamline of the simulated fluid velocity which shows, as expected, that the velocity of the fluid is higher within the fracture, thus the solute flows rapidly through the fracture see \figref{\ref{velo_concen3}}(b).

Once the vector field of the velocity is obtained, we compute the time
step $\Delta t_T$ obtained using \Equaref{\ref{local_time1}} for
$\mathbf{C}_{max} = 0.08$. The solution domain $\Omega$ is split into
sub-domain $(\Omega_i,\Delta t_i), i=0,1,2$. Here, we consider the
P\'eclet number $P_e = 3000$. 
The overlapped sub-domains are illustrated in \figref{\ref{frac}} (b).

\subsubsection{Numerical results for the case \texorpdfstring{$R(C)=0$}{Lg}}

\begin{figure}[H]
\centering
\subfloat[]{\includegraphics[height=0.18\textheight,width=0.48\textwidth]{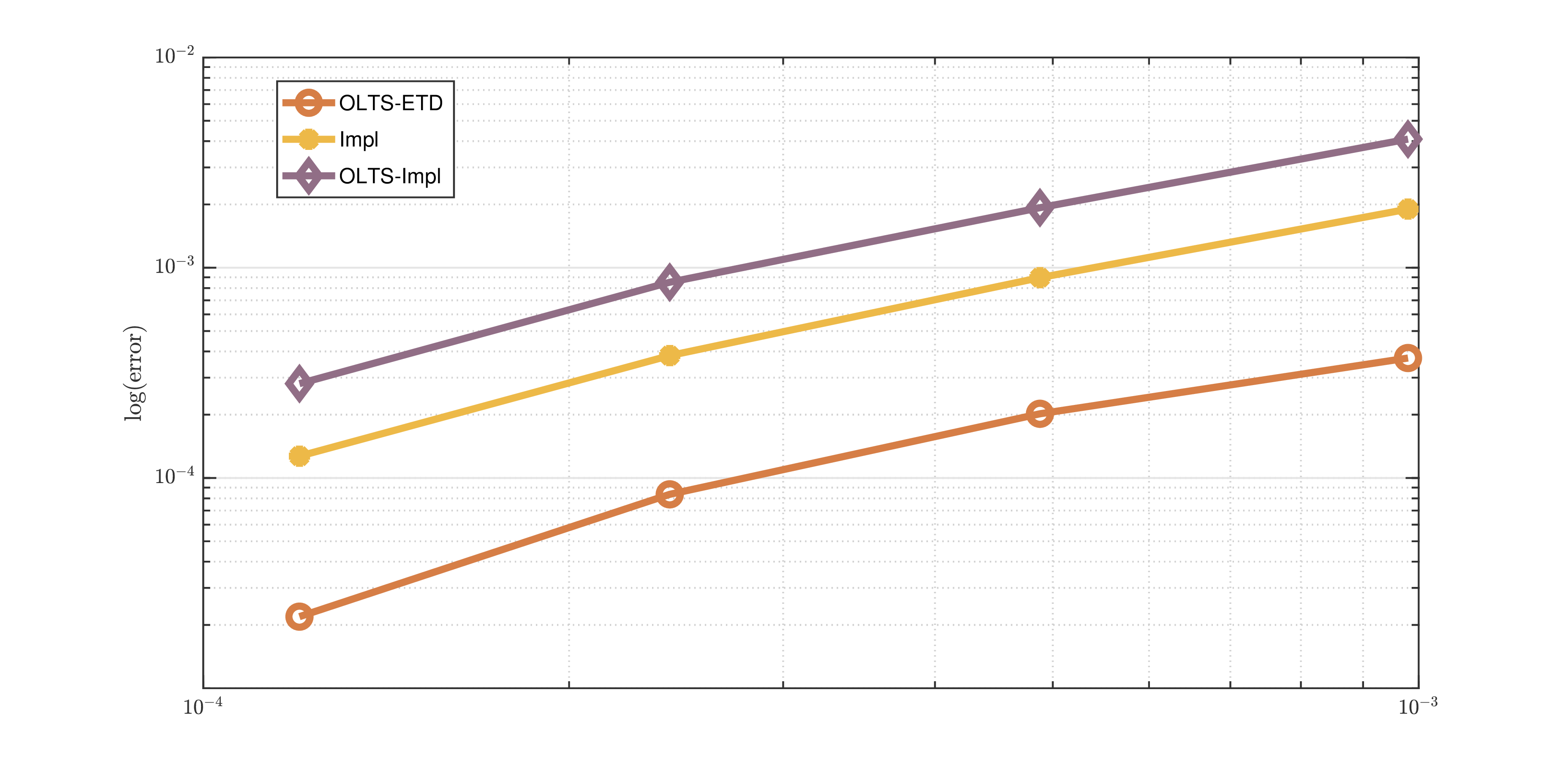}}
\subfloat[]{\includegraphics[height=0.18\textheight,width=0.48\textwidth]{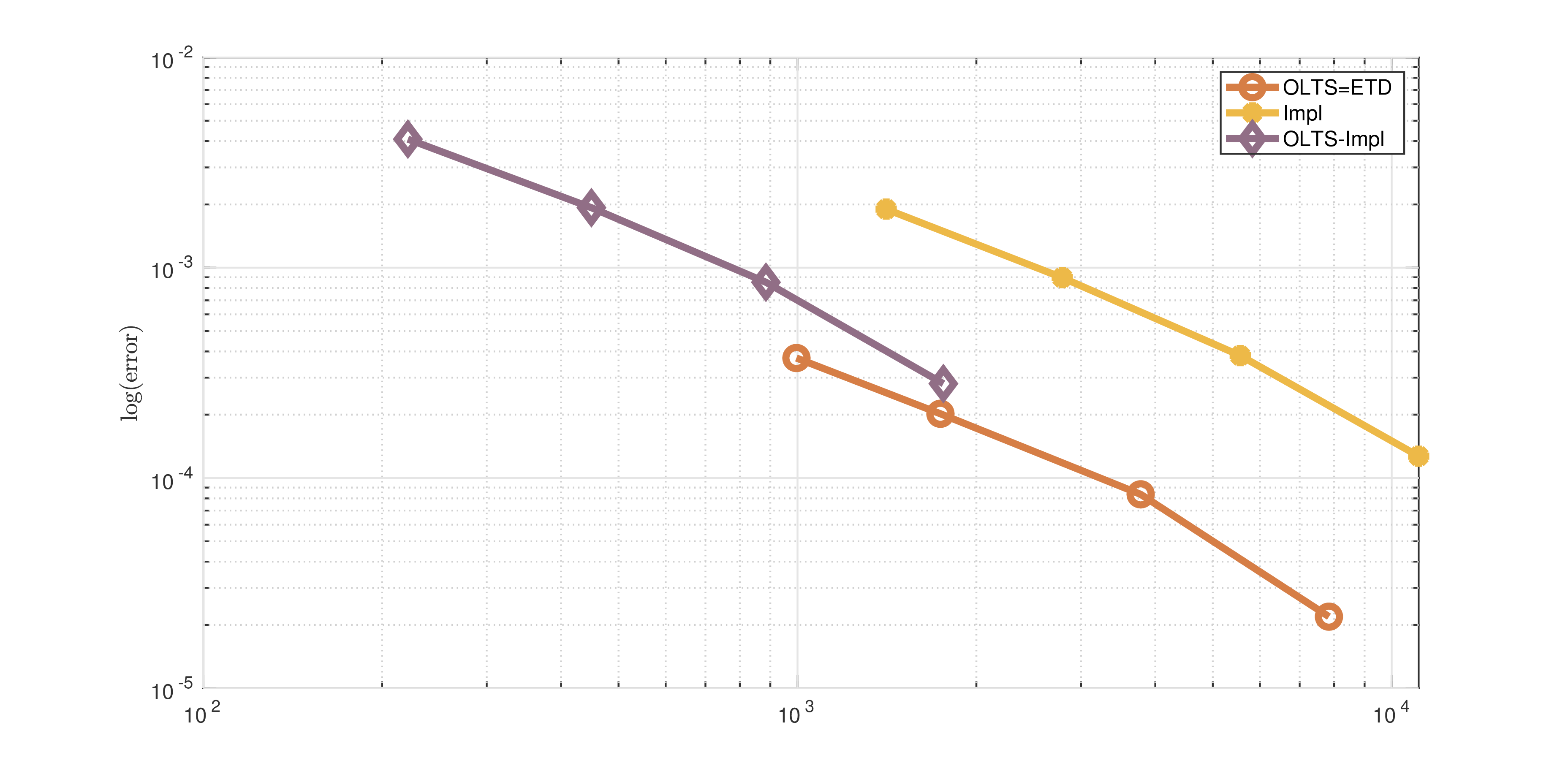}} 
\caption{Convergence and efficiency of GTS-DG and OLTS-DG schemes while solving the transport of solute through a domain with fracture without presence of source and reaction.
In (a), we plot $\log(\text{error}^r_{DG_{Impl}})$ and $\log(\text{error}^r_{DG_{OLTS-q}})$  against $\log(ht^r)$ for all $q=\text{Impl}$, $\text{ETD1}$ and $r= 0,\cdots,3$. In (b), we respectively plot $\log(\text{error}^r_{DG_{Impl}})$ and $\log(\text{error}^r_{DG_{OLTS-q}})$  against $\log(\text{CPU}^r_{DG_{q}})$ and $\log(\text{CPU}^r_{DG_{OLTS-q}})$ for all $r= 0,\cdots,3$ and $q=\text{Impl}$, $\text{ETD1}$.}
\label{compi4}
\end{figure}

For a given $r\in \{0,\cdots,4\}$, we simulate the concentration $C^r_{p}$ of the solute at the time $t=1$ using the solver $p\in\{  DG_{\text{OLTS-Impl}}, DG_{\text{OLTS-ETD1}}\}$ where the sub-domain $\Omega_i$ has the local time step $\Delta t^r_i = 2^{-r}\times \Delta t_i$ for all $i\in \{ 0,1,2\}$. Moreover, for the GTS-DG schemes ($\{DG_{\text{Impl}}\}$) we use the universal time step $ht^r$ given by
$ht^r = \min \{ \Delta t^r_i,\; i= 0,1,2  \},$
to simulate the concentration $C^r_{p}$, $r\in \{0,\cdots,4\}$ of the
solute at the time $t=1$. Throughout these simulations (i.e. for all
$r\in \{0,\cdots,4\}$), we record the computation time,
$\text{CPU}^r_{p}$, for all solver $p\in\{ DG_{Impl},
DG_{\text{OLTS}-q}\}$. To investigate the convergence, the accuracy and efficiency of solvers used here, we
assume that for a given time integrator $q\in \{Impl, ETD1  \}$, the
exact concentration of the solute at the time $t=1$ is given by
$C^4_{DG_{q}}$. We then compute the error, $\text{error}^r_{p}$, given
by  
\begin{equation*}
\text{error}^r_{DG_q} = \parallel C^4_{DG_{q}} - C^{r}_{DG_q} \parallel_{L^2(\mathcal{T})},\;\;
\text{error}^r_{DG_{OLTS-q}} = \parallel C^4_{DG_{q}} - C^{r}_{DG_{OLTS-q}} \parallel_{L^2(\mathcal{T})},
\end{equation*}
for all $r\in \lbrace 0,\cdots,3\rbrace$ and all time integrator $q\in\{ Impl, ETD1\}$.
%
%
In \figref{\ref{compi4}}(a), we plot the errors $\log(\text{error}^r_{DG_{Impl}})$ and $\log(\text{error}^r_{DG_{OLTS-q}})$  against $\log(ht^r)$ for all $q=\text{Impl}$, $\text{ETD1}$ and $r= 0,\cdots,3$. In \figref{\ref{compi4}}(b),  we plot  $\log(\text{error}^r_{DG_{q}})$ and $\log(\text{error}^r_{DG_{OLTS-q}})$  against $\log(\text{CPU}^r_{DG_{Impl}})$ and $\log(\text{CPU}^r_{DG_{OLTS-q}})$ for all $q=\text{Impl}$, $\text{ETD1}$ and $r= 0,\cdots,3$.
Note from \figref{\ref{compi4}}(a) that the errors $\text{error}^r_{p}$ decrease with the time step $ht^r$, meaning the GTS-DG and OLTS-DG schemes, considered in this section, converge. Also, note from \figref{\ref{compi4}}(a) that $\text{error}^r_{DG_{Impl}}<\text{error}^r_{DG_{OLTS-Impl}}.$
This shows that the OLTS-DG schemes is less accurate compared to GTS-DG schemes while using the same time integrator. This is expected since the GTS-DG schemes, unlike the OLTS-DG schemes, consider the finest time step, uniformly on the solution domain.
However, note from  \figref{\ref{compi4}}(b) that the computation time is reduced enough to make the OLTS-DG schemes more efficient compared to GTS-DG schemes. Moreover, \figref{\ref{compi4}}(a) shows that the accuracy improved with the accuracy of the time integrator, as $\text{error}^r_{DG_{OLTS-ETD1}}<\text{error}^r_{DG_{OLTS-Impl}}$ .

\subsubsection{Numerical results for the case \texorpdfstring{$R(C)=C-C^3$}{LG}}
\begin{figure}[H] 
\centering
\subfloat[]{\includegraphics[height=0.2\textheight,width=0.48\textwidth]{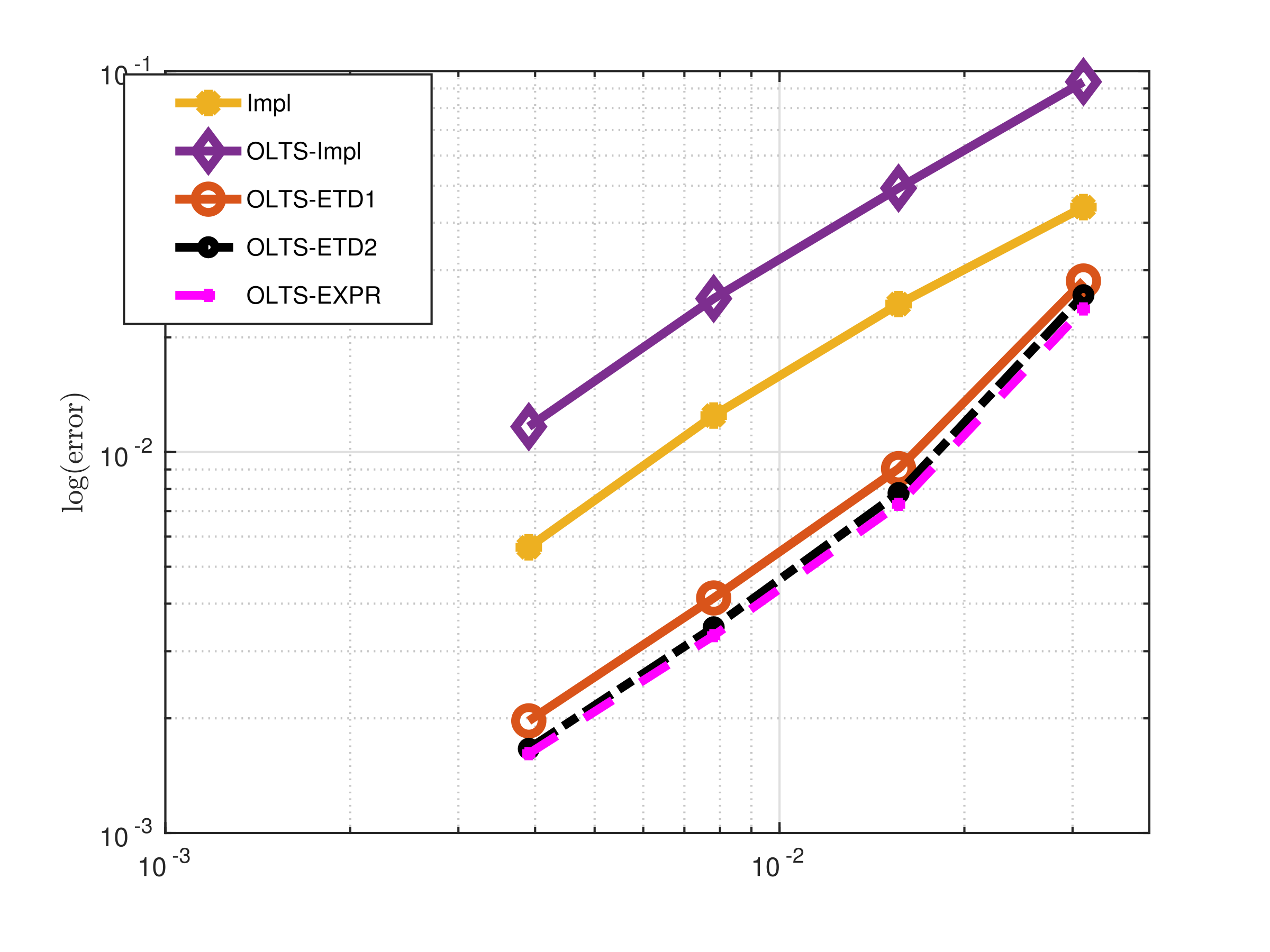}} 
\subfloat[]{\includegraphics[height=0.2\textheight,width=0.48\textwidth]{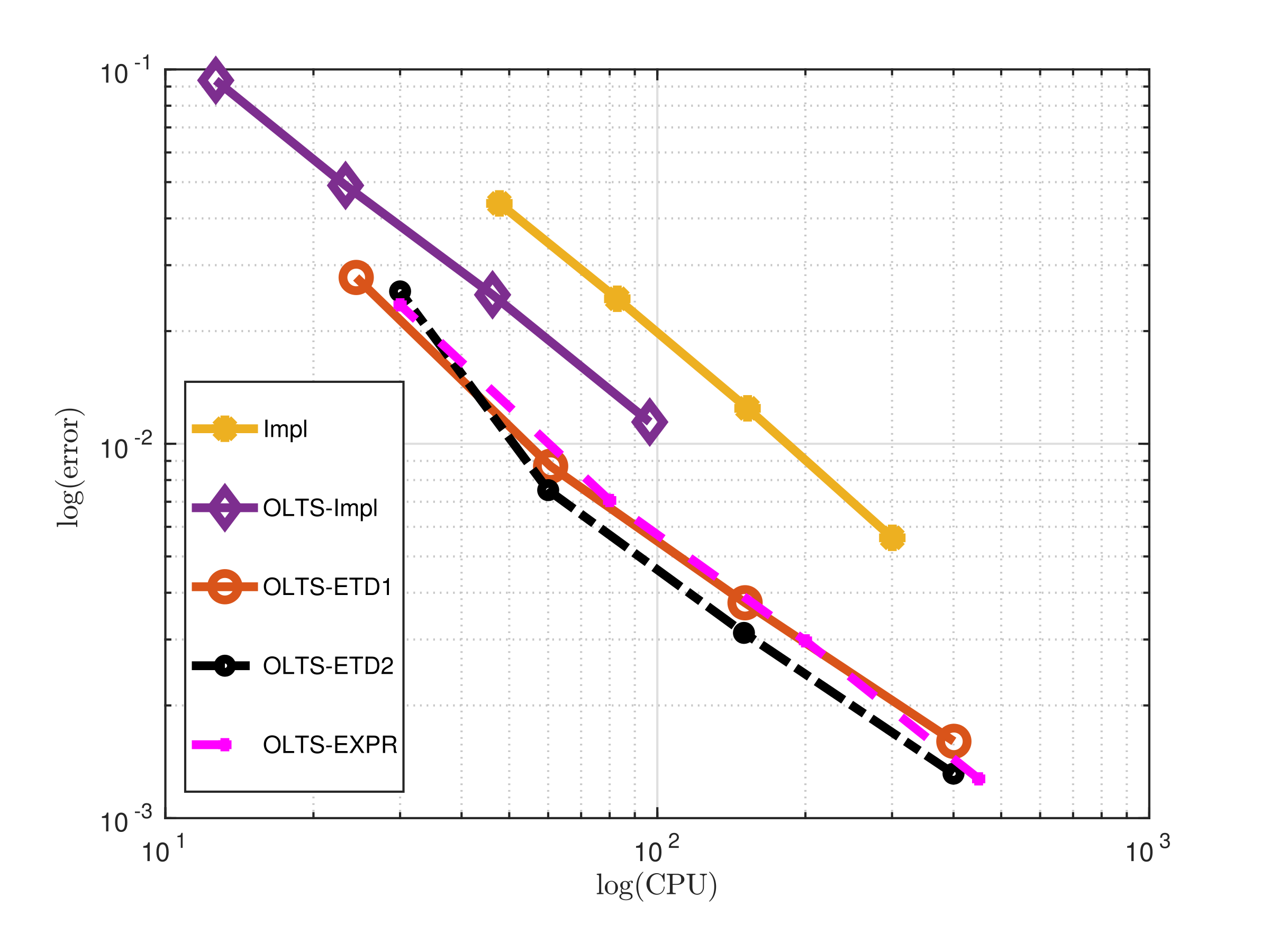}} 
\caption{Convergence and efficiency of GTS-DG and OLTS-DG schemes while solving the transport of solute through a domain with fracture with presence of reaction.
In (a), we plot $\log(\text{error}^r_{DG_{Impl}})$ and $\log(\text{error}^r_{DG_{OLTS-q}})$  against $\log(ht^r)$ for all $r= 0,\cdots,3$ and $q=\text{Impl}$, $\text{ETD1}$, $\text{ETD2}$, $\text{EXPR}$. In (b), we respectively plot $\log(\text{error}^r_{DG_{Impl}})$ and $\log(\text{error}^r_{DG_{OLTS-q}})$  against $\log(\text{CPU}^r_{DG_{Impl}})$ and $\log(\text{CPU}^r_{DG_{OLTS-q}})$ for all $r= 0,\cdots,3$ and $q\in \{\text{Impl}$, $\text{ETD1}$, $\text{ETD2}$, $\text{EXPR}\}$.}
\label{compi40}
\end{figure}

Let us consider the flow and transport of solute through a domain with
fracture with the presence of non linear reaction term, given by
$R(C)=C^3-C$. We have investigated this problem in \cite{AS17},
with the GTS-DG solvers $DG_{Impl}$ and compare the accuracy and
efficiency of the 
solvers $DG_q$ and $DG_{OLTS-q}$ with $q\in \{\text{Impl},\;
\text{ETD1},\; \text{ETD2},\; \text{EXPR}\}$. We use the previous
strategy to illustrate the results of this comparison in
\figref{\ref{compi40}}.

In \figref{\ref{compi40}}(a), we plot $\log(\text{error}^r_{DG_{Impl}})$ and
$\log(\text{error}^r_{DG_{OLTS-q}})$  against $\log(ht^r)$ for all $r=
0,\cdots,3$ and $q\in\{\text{Impl}$, $\text{ETD1}$, $\text{ETD2}$,
$\text{EXPR}\}$. In \figref{\ref{compi40}}(b), we respectively plot
$\log(\text{error}^r_{DG_{Impl}})$ and
$\log(\text{error}^r_{DG_{OLTS-q}})$  against
$\log(\text{CPU}^r_{DG_{Impl}})$ and $\log(\text{CPU}^r_{DG_{OLTS-q}})$
for all $r= 0,\cdots,3$ and $q\in\{\text{Impl}$, $\text{ETD1}$,
$\text{ETD2}$, $\text{EXPR}\}$. 

\newpage
Note from \figref{\ref{compi40}}(a) that the errors $\text{error}^r_{p}$ decrease with the time step $ht^r$, meaning the GTS-DG and OLTS-DG schemes, considered in this section, converge. Also, note from \figref{\ref{compi40}}(a) that 
\begin{equation}
\text{error}^r_{DG_{Impl}}<\text{error}^r_{DG_{OLTS-Impl}}.
\end{equation}
This shows that OLTS-DG schemes is
less accurate compared to GTS-DG schemes. This is expected since the
GTS-DG schemes, unlike the OLTS-DG schemes, consider the finest time
step, uniformly on the solution domain. However, note from
\figref{\ref{compi40}}(b) that the computation time is reduced enough to make
the OLTS-DG schemes more efficient compared to GTS-DG schemes.

\section{Conclusions}
\label{sec:conclusions}
In order to efficiently capture the localized small-scale physics of DAREs on a complex geometry, we developed here two solvers, the overlap and non overlap LTS-DG schemes, based on the domain decomposition techniques, the DG spatial discretization method and  the standard time integrators such as Impl, ETD, EXPR. The several numerical investigations lead to the following findings:
\begin{itemize}
\item When applied to Ogata and Banks problem, the numerical results of the overlap LTS-DG method show that the choice of the fast and slow components can significantly affect the accuracy of the solution. A better accuracy is obtained if the eligible sub-domains are considered in the same direction as the bulk velocity. In a high P\'eclet number regime, unlike in the low P\'eclet number regime, the size of the overlap doesn't improve the accuracy of the overlap LTS-DG method. 
\item When applied to the one dimension ETO model and the two dimension  transport of solute through a domain with fracture, the numerical results showed that the computation time is reduced enough to make the LTS-DG schemes proposed here more efficient compared to the GTS-DG schemes. These numerical results also showed that the non overlap LTS-DG method is more accurate and efficient compared to the overlap LTS-DG method. This is due to the fact that the needed information are more accurately computed in the case of the non overlap LTS method.
\end{itemize}
Note that for the LTS-DG methods proposed in this article, the same time integrator is used to advance locally the solution in time. Thus, to further improve our solvers, we will next investigate the case where different optimal time integrators will be used on different sub-domains to locally advance to solution in time.

\section*{Acknowledgments}
We would like to thank the department of applied mathematics and computer sciences of Heriot-Watt university and  Schlumberger Gould Research Centre for their support and grant of computer equipment.


%
%
%
%

\end{document}